\documentclass[10pt]{article}
\usepackage{fancyhdr}
\usepackage{extramarks}
\usepackage{amsmath}
\usepackage{amsthm}
\usepackage{amsfonts}
\usepackage{siunitx}
\usepackage{tikz}
\usepackage[plain]{algorithm}
\usepackage{algpseudocode}
\usepackage{multirow}
\usepackage{booktabs}
\usepackage{palatino}
\usepackage{subfigure}
\usepackage[colorlinks,linkcolor=black,anchorcolor=black,citecolor=black,urlcolor=blue]{hyperref}
\usepackage{amsmath,bm}
\usepackage{booktabs}
\usepackage{mathtools}
\usepackage{amssymb}
\usepackage{caption}
\usepackage{capt-of}
\usepackage{mciteplus}
\usepackage{cite}
\usepackage{mathrsfs}
\usepackage[title,titletoc,toc]{appendix}
\usepackage{xr}
\usepackage{parskip}
\usepackage{textcomp}
\usepackage[colaction]{multicol}
\usepackage[switch]{lineno}
\usepackage{lipsum}
\usepackage{etoolbox}
\usepackage{longtable}
\usepackage{array}
\usepackage{tablefootnote}
\usepackage{ragged2e}
\newcolumntype{C}[1]{>{\centering\arraybackslash}p{#1}}
\usetikzlibrary{automata,positioning}
\DeclareMathOperator{\im}{im}

\usetikzlibrary{automata,positioning}

\begin{document}

\title{HERMES: Persistent spectral graph software}

\author{Rui Wang$^1$, Rundong Zhao$^2$, Emily Ribando-Gros$^2$, Jiahui Chen$^1$, \\
Yiying Tong$^2$ \footnote{Corresponding author.		E-mail: ytong@msu.edu}, and Guo-Wei Wei$^{1,3,4}$\footnote{
Corresponding author.		E-mail: weig@msu.edu} \\
$^1$ Department of Mathematics, \\
Michigan State University, MI 48824, USA.\\
$^2$ Department of Computer Science and Engineering,\\
Michigan State University, MI 48824, USA. \\
}
\date{\today} 

\maketitle

\begin{abstract}
Persistent homology (PH) is one of the most popular tools in topological data analysis (TDA), while graph theory has had a significant impact on data science. Our earlier work introduced the persistent spectral graph (PSG) theory as a unified multiscale paradigm to encompass TDA and  geometric analysis. In PSG theory, families of persistent Laplacians (PLs) corresponding to various topological dimensions are constructed via a filtration to sample a given dataset at multiple scales. The harmonic spectra from the null spaces of PLs offer the same topological invariants, namely persistent Betti numbers, at various dimensions as those provided by PH, while the non-harmonic spectra of PLs give rise to additional geometric analysis of the shape of the data. In this work, we develop an open-source software package, called highly efficient robust multidimensional evolutionary spectra (HERMES), to enable broad applications of PSGs in  science, engineering, and technology. To ensure the reliability and robustness of HERMES, we have validated the software with simple geometric shapes and complex datasets from three-dimensional (3D) protein structures. We found that the smallest non-zero eigenvalues are very sensitive to  data abnormality.

\end{abstract}

\justifying
{\bf Availability:} HERMES is available online \url{https://weilab.math.msu.edu/HERMES/} and via GitHub \url{https://github.com/wangru25/HERMES}.

\pagenumbering{roman}
\begin{verbatim}
\end{verbatim}

{\setcounter{tocdepth}{4} \tableofcontents}
  \newpage

\setcounter{page}{1}
\renewcommand{\thepage}{{\arabic{page}}}

\section{Introduction}
As a branch of discrete mathematics, graph theory focuses on the relations among vertices  or nodes (0-simplices), edges (1-simplices), faces (2-simplices), and their high-dimensional extensions. Benefit from the capability of graph formulations that encode inter-dependencies among constituents of versatile data into simple representations, graph theory has been regarded as the mathematical scaffold in the study of various complex systems in biology, material science,  physical infrastructure, and network science.  However,  traditional graphs only represent the pairwise relationships between entries. Therefore, hypergraphs, a generalization of graphs that describe the multi-way relationships of mathematical structures have been developed to capture the high-level complexity of data \cite{aksoy2020hypernetwork,bressan2016embedded}.  Mathematically,  graphs and hypergraphs are intrinsically related to the simplicial complexes, which have broader use in computational topology. Moreover, many other areas such as algebra, group theory, knot theory, spectral graph theory (SGT), algebraic topology (AT), and combinatorics are closely related to graph theory. Among them, the applications of SGT have been driven by various real-life problems in chemistry, physics, and life science in the past few decades \cite{nguyen2019agl, spielman2007spectral}.

In its early days, the spectral graph theory studied the properties of a graph by its graph Laplacian matrix and adjacent matrix. Later on, developments in the spectral graph theory involve some geometric flavor. The explicit constructions of expander graphs rely on studying the eigenvalues and isoperimetric properties of graphs. The discrete analog of Cheeger's inequality for graphs in Riemannian geometry is related to the study of manifolds \cite{cheeger1969lower}. Specifically, an eigenvalue of the Laplacian of a manifold is related to the isoperimetric constant of the manifold, which motivates the study of graphs by employing manifolds. Benefiting from the increasingly rich connections with differential geometry, the spectral graph theory has entered a new era \cite{chung1997spectral}. One of the critical developments is the Laplacian on a compact Riemannian manifold in the context of the de Rham-Hodge theory \cite{kamber1987rham,zhao2020rham}. The harmonic part of the Hodge Laplacian spectrum contains the topological information, whereas the non-harmonic part of the Hodge Laplacian spectrum offers additional geometric information for shape analysis \cite{chen2019evolutionary}. Indeed, the connectivity of a graph/topological space can be revealed from  topological invariants. It is well-known that the number of the eigenvalues in the harmonic spectra of $q$th-order persistent Laplacian represents the dimension of persistent $q$-cohomology of a graph \cite{hernandez2019higher,wang2020persistent}, which builds the connection between spectral graph theory and algebraic topology.

Homology and cohomology are key concepts in the algebraic topology, which were developed to analyze and classify manifolds according to their cycles. The traditional homology is genuinely metric-independent, indicating that the geometric information is barely considered \cite{kaczynski2006computational}. Therefore, for practical computation, a new branch of algebraic topology named persistent homology (PH) \cite{zomorodian2005computing,carlsson2009zigzag, edelsbrunner2008persistent} is implemented to create a sequence of topological spaces characterized by a filtration parameter, such as the radius of a ball or the level set of a real-valued function.  As the most important realization of the topological data analysis (TDA) \cite{de2007coverage, dey2014computing,bubenik2007statistical}, topological persistence has had great success in computational biology \cite{xia2014persistent,cang2017topologynet}. For instance, the superior performance of using PH features of protein-drug complexes in the free energy prediction and ranking at D3R Grand Challenges (\url{https://drugdesigndata.org/about/grand-challenge}), a worldwide competition series in computer-aided drug design  \cite{nguyen2019mathematical}, was a remarkable success for TDA.  Additionally, a weighted persistent homology is proposed as a unified paradigm for the analysis of the biomolecular data system \cite{meng2020weighted}.

Recently, we have introduced persistent spectral graph (PSG) theory to bridge persistent homology and spectral graph theory \cite{wang2019persistent,wang2020persistent}. The PSG theory extends the persistence notion or multiscale analysis to algebraic graph theory. A family of spectral graphs induced by a filtration overcome the difficulty of using traditional spectral graph theory in analyzing graph structures with a single geometry, giving rise to persistent spectral analysis (PSA). Additionally, the evolution of the null space dimension of the persistent Laplacian (PL) over the filtration offers the topological persistence. Therefore, PSG theory provides simultaneous TDA and PSA. Specifically, by varying a filtration parameter, a series of $q$th-order persistent Laplacians  (or $q$-persistent Laplacian) provide  persistent spectra. Notably, the persistent harmonic  spectra of $0$-eigenvalues span the null space of the $q$-th order persistent Laplacian and fully recover the persistent $q$-th Betti numbers or persistent barcodes \cite{carlsson2005persistence} of the associated persistent homology. Specifically, the number of $0$-eigenvalues of $q$th-order persistent Laplacian reveals the number of $q$-cocycles for a given point-cloud dataset. Moreover, the additional geometric shape information of the data will be unveiled in the non-harmonic spectra. Recently, the theoretical properties and algorithms of PSGs have been further studied  \cite{memoli2020persistent} and the application of PSG methods to drug discovery  has been reported \cite{meng2020persistent} . The de Rham-Hodge theory counterpart, called evolutionary de Rham-Hodge theory, has also been formulated \cite{chen2019evolutionary}.

Currently, many  open-source packages have been developed for the applications of persistent homology, including Ripser \cite{bauer2017ripser}, Dionysus \cite{morozov2012dionysus}, Gudhi \cite{gudhi2015gudhi}, Perseus \cite{mischaikow2013morse}, DIPHA \cite{bauer2014dipha}, Javaplex \cite{adams2014javaplex}, CliqueTop \cite{giusti2015clique}, DioDe  \cite{morozov2020diode}, \href{https://bitbucket.org/grey_narn/hera/src/master/}{Hera},  \href{http://gregoryhenselman.org/eirene/index.html}{Eirene}, and ``TDA'' package in R \cite{fasy2019package}. These packages are able to construct a family of complexes with the point clouds data as input and calculate its corresponding Betti numbers, which are equivalent to  the harmonic spectra of the persistent Laplacian. However, there is no software package for simultaneous TDA and PSA. While we developed the theoretical part of the persistent spectral graph in 2019, we have not constructed an efficient and robust software yet.

The objective of the present work is to provide the first open-source package, called highly efficient robust multidimensional evolutionary spectra  (HERMES), for evaluating both the harmonic and non-harmonic spectra of persistent Laplacians. In the present release, we consider an implementation in alpha complexes \cite{edelsbrunner2010alpha}. To verify the reliability of HERMES, 15 complicated 3D structures of proteins are used to calculate the spectra of $q$th-order persistent Laplacians for $q=0,1,2$. Moreover, as a validation, the persistent harmonic spectra generated by HERMES are compared with those obtained from Gudhi and DioDe.

\section{Method}
As a powerful and versatile data representation that encodes inter-dependencies among constituents, graph theory has widely spread applications in various fields such as molecular sciences, engineering, physics, biology, algebra, topology, and  combinatorics. In this section, we first briefly review the concepts of simplex, simplicial complex, chain complex, Delaunay complex, and alpha complex in topology, which can be regarded as generalizations of graph into its higher-dimensional topological counterparts. Then, we review the $q$th-order Laplacian for simplicial complexes, which is a generalization of the graph Laplacian  in graph theory. The topological and geometric information of a single configuration can be evaluated from the spectra of the $q$th-order Laplacian. Moreover, built upon these concepts, we will discuss persistent spectral graph \cite{wang2019persistent,wang2020persistent}
for the analysis of topological invariants and geometric measurements of high-dimensional datasets. Instead of analyzing the spectra for only one configuration, the persistent spectral graphs can analyze a series of topological and geometric changes, which enriches the set of available representations for high-dimensional datasets.

\subsection{Topological concepts}
In this section, we give a concise review of simplex, simplicial complex, and chain complex to provide essential background for   persistent spectral graphs. More details can be found in the literature \cite{zomorodian2005computing,carlsson2009zigzag, edelsbrunner2008persistent}.

\justifying
{\bf Simplex.} A $q$-simplex denoted as $\sigma_q$ is the convex hull of $q + 1$ affinely independent points in $\mathbb{R}^n$, having dimension $\text{dim}(\sigma_q) = q$. For example, a vertex is a $0$-simplex, an edge is a $1$-simplex, a triangle is a $2$-simplex, and a tetrahedron is a $3$-simplex. We call the convex hull of each non-empty subset of $q+1$ points a face of $\sigma_q$, and each of its corner points is also called one of its vertices.

\justifying
{\bf Simplicial complex.} A set of simplices is a simplicial complex denoted as $K$ if the following conditions are satisfied:
\begin{itemize}
    \item[(1)] If all faces of any simplex in $K$ are also in $K$, and
    \item[(2)] The non-empty intersection of any two simplices in $K$ is a common face of the two simplices.
\end{itemize}
The dimension of simplicial complex $K$ is defined as $\text{dim}(K) = \max\{\text{dim}{\sigma_q}: \sigma_q \in K\}$.

\justifying
{\bf Chain complex.} A $q$-chain is a formal sum of $q$-simplices in simplicial complex $K$ with $\mathbb{Z}_2$ coefficients. The set of all $q$-chains has a basis which the set of $q$-simplices in $K$, thus forming a finitely generated free abelian group denoted as $C_q(K)$. The boundary operator is a group homomorphism defined by $\partial_q: C_q(K) \to C_{q-1}(K)$ to relate the chain groups. More specifically, denoting $q$-simplex as $\sigma_q = [v_0, v_1, \cdots, v_q]$ by its vertices $v_i$, the boundary operator is defined through its action on the basis,
\begin{equation}
    \partial_q \sigma_q = \sum_{i=0}^{q}(-1)^i\sigma^{i}_{q-1}.
\end{equation}
Here, $\sigma^{i}_{q-1} = [v_0, \cdots, \hat{v_i},\cdots,v_q]$ is the $(q\!-\!1)$-simplex with $v_i$ omitted. The following sequence of chain groups connected by boundary operators is a \emph{chain complex} (defined as a set of abelian groups connected by homomorphisms such that the composite of any two consecutive homomorphisms is zero, $\partial_q\partial_{q+1}=0.$)
\[
\cdots \stackrel{\partial_{q+2}}\longrightarrow C_{q+1}(K) \stackrel{\partial_{q+1}}\longrightarrow C_{q}(K) \stackrel{\partial_{q}}\longrightarrow C_{q-1}(K)\stackrel{\partial_{q-1}} \longrightarrow \cdots
\]

\subsection{Combinatorial Laplacians}\label{sec:qLaplacian}
Combinatorial Laplacians offer both spectral analysis and topological analysis \cite{hernandez2019higher}. One central role played by the chain complex associated with a simplicial complex is to define its $q$-th homology group ($H_q=\ker \partial_q/\im \partial_{q+1}$), which is a topological invariant of the simplicial complex. The dimension of $H_q$ is denoted by $\beta_q=\dim H_q$, the $q$-th Betti number, which, roughly speaking, measures the number of $q$-dimensional holes in the simplicial complex, or the geometric object tessellated into the simplicial complex.

A dual chain complex can be defined on any chain complex through the adjoint operator of $\partial_q$ defined on the dual spaces $C^q(K)=C^{\ast}_q(K).$ The $q$-coboundary operator $\partial_q^{\ast}: C^{q-1}(K) \to C^q(K)$ is defined as:
\begin{equation}
    \partial^\ast \omega^{q-1}(c_{q}) \equiv \omega^{q-1}(\partial c_q),
\end{equation}
where $\omega^{q-1}\in C^{q-1}(K)$ is a \emph{$(q\!-\!1)$-cochain}, which is a homomorphism mapping a chain to the coefficient group, and $c_q\in C_q(K)$ is a $q$-chain. The homology of the dual chain complex is often called \emph{cohomology}.

If we denote by $\mathcal{B}_q$ the matrix representation of a $q$-boundary operator with respect to the standard basis for $C_q(K)$ and $C_{q-1}(K)$, the number of rows and the number of columns in $\mathcal{B}_q$ correspond to the number of $(q-1)$-simplices and that of $q$-simplices in $K$, respectively. Moreover, the matrix representation of $q$-coboundary operator is denoted $\mathcal{B}_q^{T}$.

In de Rham-Hodge theory, homology and cohomology are often studied through their correspondences to the $q$-combinatorial Laplacian operator, defined as the linear operator $\Delta_q: C^q(K) \to C^q(K)$ as follows,
\begin{equation}
    \Delta_q := \partial_{q+1}\partial_{q+1}^{\ast} + \partial_q^{\ast}\partial_q,
\end{equation}
where the isomorphism $C^q(K)\cong C_q(K)$ is assumed, where each $q$-simplex is mapped to its own dual, i.e., the isomorphism keeps the coefficients of chains and cochains in the standard simplicial basis. Correspondingly, the matrix representation of $\Delta_q$ is the  $q$th-order Laplacian, which is denoted $\mathcal{L}_q(K)$,
\begin{equation}
    \mathcal{L}_q(K) =  \mathcal{B}_{q+1}\mathcal{B}_{q+1}^{T} + \mathcal{B}_q^{T} \mathcal{B}_q.
\end{equation}
Assume the number of $q$-simplices exist in $K$ to be $N_q$, then $\mathcal{L}_q(K)$ is an $N_q\!\!\times\!\!N_q$-matrix. Since the $q$th-order Laplacian $\mathcal{L}_q(K)$ is symmetric and positive semi-definite, its spectrum consists of only real and non-negative eigenvalues. We denote the spectrum of $\mathcal{L}_q(K)$ as
\[
\text{Spec}(\mathcal{L}_q(K)) = \{\lambda_{1,q}, \lambda_{2,q}, \cdots, \lambda_{N_q,q}\}.
\]
The multiplicity of zero in the spectrum (also called the harmonic spectrum) reveals the topological information $\beta_q$, whereas the non-harmonic spectrum encodes further geometric information. The correspondence between the multiplicity of zero spectra of $\mathcal{L}_q(K)$ and the $q$th Betti number defined in the homology is an important result in de Rham-Hodge theory,
\cite{kamber1987rham,zhao2020rham,chen2019evolutionary}
\begin{equation}
    \beta_q = \dim \ker \partial_q-\dim \im \partial_{q+1}= \dim \ker \mathcal{L}_q(K) = \# \text{0 eigenvalues of } \mathcal{L}_q(K).
\end{equation}
Intuitively, $\beta_0$ represents the number of connected components in $K$, $\beta_1$ reveals the number of 1D noncontractible loops or circles in $K$, and $\beta_2$ shows the number of 2D voids or cavities in $K$.

\subsection{Persistent spectral graphs}
Both topological and geometric information can be derived from analyzing the spectra of $q$th-order Laplacian. However, the information is restricted to those pieces contained in the connectivity of the simplicial complex. A single simplicial complex produces insufficient information for practical problems such as feature extraction for machine learning analysis. To enrich the spectral information, persistent spectral graph (PSG) is proposed by creating a sequence of simplicial complexes induced by varying a filtration parameter, which is inspired by persistent homology as well as our earlier multiscale graph Laplacians \cite{xia2015multiscale}.

First, we consider a filtration of simplicial complex $K$ which is a nested sequence of subcomplexes $(K_t)_{t=0}^m$ of the final complex $K$:
\begin{equation}
    \emptyset = K_0 \subseteq K_1 \subseteq K_2 \subseteq \cdots \subseteq K_m = K.
\end{equation}
For each subcomplex $K_t$, we denote its corresponding chain group to be $C_q(K_t)$, and the $q$-boundary operator will be denoted by $\partial_q^t: C_q(K_t) \to C_{q-1}(K_t)$. As conventionally done, we define $C_q(K_t)$ for $q<0$ as the zero group $\{0\}$ and $\partial_q^t$ as a zero map. \footnote{We define the boundary matrix $\mathcal{B}_0^t$ for the boundary operator $\partial_0^t$ as a zero matrix. The number of columns of $\mathcal{B}_0^t$ is the number of $0$-simplices in $K_t$, the number of rows will be $1$.} If $0<q \le \dim K_t$, then
\begin{equation}
    \partial_q^t(\sigma_q) = \sum_{i}^q(-1)^i \sigma^i_{q-1}, \quad \forall \sigma_q \in K_t,
\end{equation}
with $\sigma_q = [v_0, \cdots, v_q]$ being any $q$-simplex, and $\sigma^{i}_{q-1} = [v_0, \cdots, \hat{v_i} ,\cdots,v_q]$ being the $(q\!-\!1)$-simplex constructed by removing $v_i$ . The adjoint operator of $\partial_q^t$ is the coboundary operator $\partial_q^{t^{\ast}}: C^{q-1}(K_t) \to C^q(K_t),$  which can be regarded as a map from $C_{q-1}(K_t)$ to $C_q(K_t)$ through the isomorphisms $C^q(K_t)\cong C_q(K_t)$ between cochain groups and chain groups.

Similar to the persistent homology, a sequence of chain complexes can be defined as below:
\begin{equation}
    \left.\begin{array}{cccccccccccccc}
        \cdots & C_{q+1}^1 &
        \xrightleftharpoons[\partial_{q+1}^{1^\ast}]{\partial_{q+1}^1} & C_q^1 &
        \xrightleftharpoons[\partial_q^{1^\ast}]{\partial_q^1} & \cdots & \xrightleftharpoons[\partial_3^{1^\ast}]{\partial_3^1} & C_2^1 & \xrightleftharpoons[\partial_2^{1^\ast}]{\partial_2^1} & C_1^1 & \xrightleftharpoons[\partial_1^{1^\ast}]{\partial_1^1} & C_0^1 & \xrightleftharpoons[\partial_0^{1^\ast}]{\partial_0^1} & C_{-1}^1=\{0\} \\
        & \rotatebox{-90}{$\subseteq$} &  & \rotatebox{-90}{$\subseteq$} &  &  &  & \rotatebox{-90}{$\subseteq$} &  & \rotatebox{-90}{$\subseteq$} &  & \rotatebox{-90}{$\subseteq$} &  &  \\
        \cdots & C_{q+1}^2 &
        \xrightleftharpoons[\partial_{q+1}^{2^\ast}]{\partial_{q+1}^2} & C_q^2 &
        \xrightleftharpoons[\partial_q^{2^\ast}]{\partial_q^2} & \cdots &
        \xrightleftharpoons[\partial_3^{2^\ast}]{\partial_3^2} & C_2^2 & \xrightleftharpoons[\partial_2^{2^\ast}]{\partial_2^2} & C_1^2 & \xrightleftharpoons[\partial_1^{2^\ast}]{\partial_1^2} & C_0^2 & \xrightleftharpoons[\partial_0^{2^\ast}]{\partial_0^2} & C_{-1}^2=\{0\} \\
        & \vdots & & \vdots &&&&\vdots &&\vdots && \vdots && \vdots\\
        \\
        & \rotatebox{-90}{$\subseteq$} &  & \rotatebox{-90}{$\subseteq$} &  &  &  & \rotatebox{-90}{$\subseteq$} &  & \rotatebox{-90}{$
        \subseteq$} &  & \rotatebox{-90}{$\subseteq$} &  &  \\
        \cdots & C_{q+1}^m &
        \xrightleftharpoons[\partial_{q+1}^{m^\ast}]{\partial_{q+1}^m} & C_q^m &
        \xrightleftharpoons[\partial_q^{m^\ast}]{\partial_q^m} & \cdots &
        \xrightleftharpoons[\partial_3^{m^\ast}]{\partial_3^m} & C_2^m & \xrightleftharpoons[\partial_2^{m^\ast}]{\partial_2^m} & C_1^m & \xrightleftharpoons[\partial_1^{m^\ast}]{\partial_1^m} & C_0^m & \xrightleftharpoons[\partial_0^{m^\ast}]{\partial_0^m} & C_{-1}^m=\{0\}
    \end{array}\right.
\end{equation}
For simplicity, we use $C_q^t$ to denote the chain group $C_q(K_t)$.

Next, we introduce persistence to the Laplacian spectra. We define the subset of $C_q^{t+p}$ whose boundary is in $C_{q-1}^t$ as $\mathbb{C}_q^{t,p}$, assuming the natural inclusion map from $C_{q-1}^t$ to $C_{q-1}^{t+p}$,
\begin{equation}
    \mathbb{C}_q^{t,p} \coloneqq \{ \beta \in C_q^{t+p} \ | \ \partial_q^{t+p}(\beta) \in C_{q-1}^{t}\}.
\end{equation}
On this subset, one may define the $p$-persistent $q$-boundary operator denoted by $\eth_q^{t,p} : \mathbb{C}_q^{t,p} \to  C_{q-1}^{t}$. Its corresponding adjoint operator is $(\eth_q^{t,p})^{\ast} : C_{q-1}^{t}  \to  \mathbb{C}_q^{t,p}$, again through the identification of  cochains with chains. We then define the  $q$-order $p$-persistent Laplacian operator $\Delta_q^{t,p}: C_q^t \to C_q^t$ associated with the filtration as
\begin{equation}
    \Delta_q^{t,p} = \eth_{q+1}^{t,p} \left(\eth_{q+1}^{t,p}\right)^\ast + \partial_q^{t^\ast} \partial_q^t.
\end{equation}
The matrix representation of $\Delta_q^{t,p}$ in the simplicial basis is
\begin{equation}
    \mathcal{L}_q^{t,p} = \mathcal{B}_{q+1}^{t,p} (\mathcal{B}_{q+1}^{t,p})^T + (\mathcal{B}_{q}^t)^T \mathcal{B}_{q}^t,
\end{equation}
where $\mathcal{B}_{q+1}^{t,p}$ is the matrix representation of $\eth_{q+1}^{t,p}.$

We denote the spectrum of $\mathcal{L}_q^{t,p}$ as
\[
\text{Spec}(\mathcal{L}_q^{t,p}) = \{\lambda_{1,q}^{t,p}, \lambda_{2,q}^{t,p}, \cdots, \lambda_{N_q^t,q}^{t,p}  \},
\]
where $N_q^t=\dim C_q^{t}$ is the number of $q$-simplices in $K_t$, and the eigenvalues are listed in the ascending order. Thus, the smallest non-zero eigenvalue of $\mathcal{L}_q^{t,p}$ is denoted as $\lambda_{2,q}^{t,p}$. We may recognize the multiplicity of zero in the spectrum of $\mathcal{L}_q^{t,p}$ as the $q$th order $p$-persistent  Betti number $\beta_q^{t,p}$, which counts the number of (independent) $q$-dimensional holes in $K_t$ that still exists in $K_{t+p}$. The relation can be observed in
\begin{equation}
    \beta_q^{t,p} = \dim \ker \partial_q^t - \dim \im \eth_{q+1}^{t,p} = \dim \ker \mathcal{L}_q^{t,p} = \# \text{0 eigenvalues of } \mathcal{L}_q^{t,p}.
\end{equation}
In this paper, we focus on the $0,1,2$th-order persistent Laplacians, which depict the relations among vertices, edges, triangles, and tetrahedra, as we target 3D real-world applications.

For instance, given a set of vertices $V =\{v_0, v_1, \cdots, v_{N_0-1}\}$ , $N_0$ embedded in $\mathbb{R}^3$, we consider a nested family of simplicial complexes that may be created for a positive real number $\alpha$. Denoting the simplicial complex generated for $\alpha$ by $K_{\alpha}$, the traditional $q$th-order Laplacian is just a special case of $q$th-order $0$-persistent  Laplacian at $K_{\alpha}$
\begin{equation}
        \mathcal{L}_q^{\alpha,0} = \mathcal{B}_{q+1}^{\alpha,0}(\mathcal{B}_{q+1}^{\alpha,0})^{T} + (\mathcal{B}_q^{\alpha})^T \mathcal{B}_q^{\alpha}.
\end{equation}
The spectrum of $\mathcal{L}_q^{\alpha,0}$ is simply associated with a snapshot of the filtration,
\begin{equation}
    \text{Spec}(\mathcal{L}_q^{\alpha,0})  = \{\lambda_{1,q}^{\alpha,0}, \lambda_{2,q}^{\alpha,0}, \cdots, \lambda_{N_q^{\alpha},q}^{\alpha,0}  \}.
\end{equation}
Correspondingly, the $q$-th $0$-persistent  Betti number $\beta_q^{\alpha,0}=\beta_q^{\alpha}$. In addition to the traditional homology information, and persistent homology information, our proposed persistent spectral graph theory, through the nonzero eigenvalues in the spectrum of the persistent Laplacian operator, provide richer spatial information induced by varying the filtration parameters. Thus it provides a powerful tool to encode high-dimensional datasets into various topological and geometric features in a coherent fashion.\footnote{In this work, we use notations $\mathbb{C}_q^{t,p}, \eth_q^{t,p}, \Delta_q^{t,p}, \mathcal{L}_q^{t,p}$, and $\beta_q^{t,p}$ instead of $\mathbb{C}_q^{t+p}, \eth_q^{t+p}, \Delta_q^{t+p}, \mathcal{L}_q^{t+p}$, and $\beta_q^{t+p}$ used in Ref. \cite{wang2020persistent}.}

\subsection{Delaunay triangulation and alpha shape}
In this section, we provide the details on a practical construction of filtration for persistent spectral graph theory based on the alpha complex. The alpha complex can be regarded as a simplicial complex, which is a homotopy equivalent to the nerve of balls around data points. Its geometric realization built as the union of convex hulls of points in each simplex is called the alpha shape. The alpha shape was first proposed in 1983, which defined the shape associated with a finite set of points in the plane controlled by one parameter~\cite{edelsbrunner2010alpha}.

In the following, we first describe how to construct the alpha shape, and then provide some necessary concepts for the implementation of the alpha complex in PSG theory. Let $P$ be a finite set of points in $q$D Euclidean space $\mathbb{R}^{q}$ ($q=2$ or $3$ in most applications), and $\alpha$ be a positive real number. Denote an open ball with radius $\alpha$ as an alpha ball ($\alpha$-ball). We say that an $\alpha$-ball is empty if it contains no point of $P$, and the alpha hull ($\alpha$-hull) of $P$ is the set of points that do not belong to any empty $\alpha$-ball. For any subset $T \subseteq P$ with size $|T| = k+1, 0\le k \le q$, the geometric realization of $k$-simplex $\sigma_{T}$ is the convex hull of $T$. We say that a $k$-simplex $\sigma_{T}$ is $\alpha$-exposed if there exists an  empty $\alpha$-ball ${\bf b}$ such that $T = \partial {\bf b} \cap P $ for $0\le k \le q-1$. Denoting the collection of $\alpha$-exposed $k$-simplices as $F_{k,\alpha}$ for $0\le k \le q-1$, the alpha shape ($\alpha$-shape) of $P$ is the polytope whose boundary consists of the $k$-simplices in $F_{k,\alpha}$. The alpha complex is just the simplicial complex that is the collection of the simplices in the alpha shape.

There are two structures that are closely related to the alpha shape and helpful in efficient implementation of alpha shape and alpha complex. One is the Voronoi diagram \cite{voronoi1908nouvelles} and the other is its dual structure, the Delaunay tessellation \cite{delaunay1934sphere}. The latter is the alpha complex for sufficiently large $\alpha$, e.g., when $\alpha$ is greater than the diameter of $P$. Thus, the Delaunay tessellation is the final complete simplicial complex in the filtration that we use.

For a given set of points $P = \{p_1, p_2, \cdots, p_n\} \subseteq \mathbb{R}^q$, the Voronoi cell $V_i$ of a point $p_i \in P$ contains all of the points for which $p_i$ is the closest among all the points in $P$,
\begin{equation}
    V_i = \{x\in \mathbb{R}^q \ | \quad \|x-p_i\| \le \|x-p_j\|,\quad  \forall p_j \in P \}.
\end{equation}
The Voronoi diagram of $P$ is the set of Voronoi cells, which is defined as
\begin{equation}
    \text{Vor} P = \{V_i \ | \ \forall i \in \{1,2,\cdots,|P|\}\}.
\end{equation}
The Delaunay tessellation for a given set $P$ in general position (i.e., no $q+1$ ponits are in a $(q\!-\!1)$-D linear subspace, and no $q+2$ points share the same circumsphere) is the dual simplicial complex to the Voronoi diagrams. For instance, a Delaunay tessellation for a given set $P$ in 2D is a triangulation $\text{DT}(P)$ such that no point in $P$ is inside the circumcircle of any triangle in $\text{DT}(P)$ \cite{aurenhammer2013voronoi,may2018multivariate}. A formal way to define the Delaunay tessellation is to use the nerve of the collection of Voronoi cells ($\text{Nrv}(\text{Vor} P)$), which can be expressed as
\begin{equation}
    \text{DT}(P) = \text{Nrv}(\text{Vor} P) = \{ J\subseteq \{1,2,...,|P|\} \;\; | \;\;  \bigcap_{i\in J} V_i \neq \emptyset  \},
\end{equation}
under the condition that the points in $P$ are general position. Note that, in practice, a set of points that are not in  general position can be symbolically perturbed to general position.

\begin{figure}[t]
    \centering
	\includegraphics[width=1\textwidth]{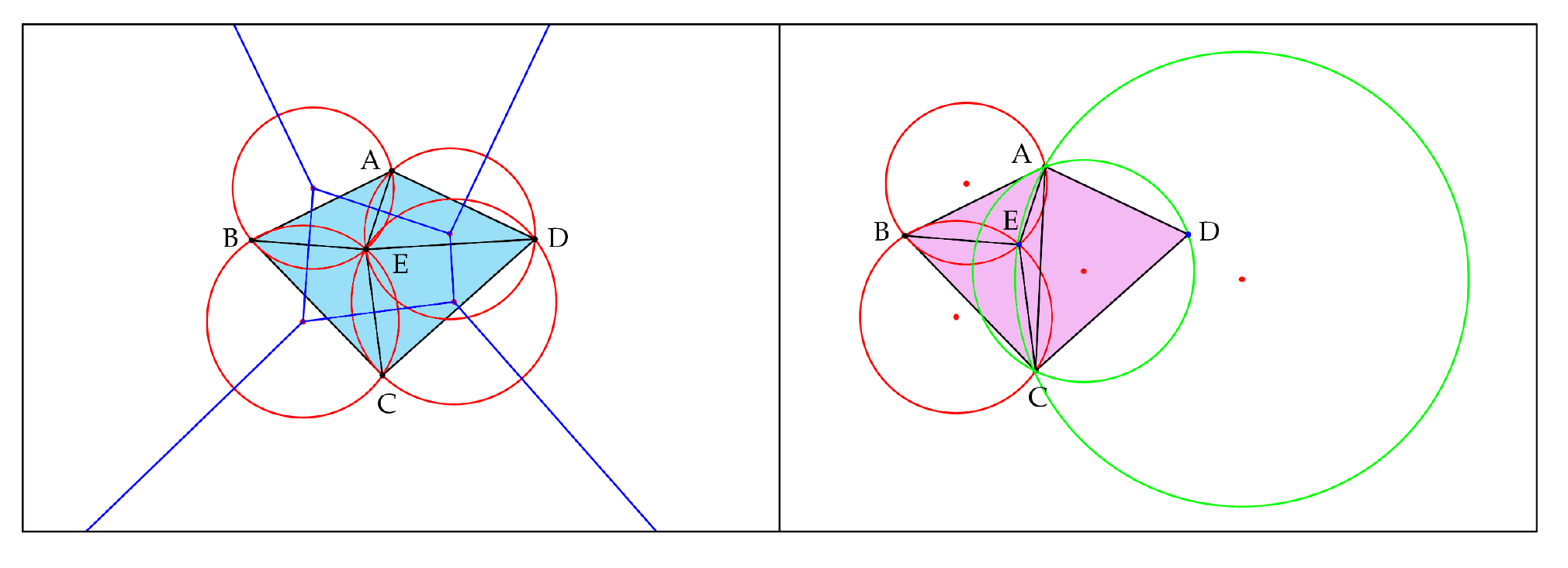}
    \caption{Illustration of Voronoi diagram, Delaunay triangulation, and Non-Delaunay triangulation. {\bf Left chart:} The Voronoi diagram and its dual Delaunay triangulation. The points set is $P$ = \{A,B,C,D,E\} and the Delaunay is defined as $\text{DT}(P)$. The blue lines tessellate the plane into Voronoi cells. The red circle are the circumcircles of triangles in $\text{DT}(P)$. {\bf Right chart:} A Non-Delaunay triangulation. Vertices E and D are in the green circumcircles, implying the right chart is an example of Non-Delaunay triangulation. }
    \label{fig:Delaunay}
\end{figure}

Next, we introduce the mathematical description of the construction of alpha complex through the union of balls centered at points in $P$, which is essentially a van der Waals surface for atoms positioned at $P$ with the same radius $\alpha$. For a given set of points $P = \{p_1, p_2, \cdots, p_n\}$ in $\mathbb{R}^q$ and a positive  real number $\alpha$, we can denote the closed ball centered at $p_i$ as $B_i(\alpha) = p_i + \alpha\mathbb{B}^q,$ where $\mathbb{B}^q$ is a $q$D unit ball around the origin. The union of these balls can be expressed as
\begin{equation}
    U(\alpha) = \{x\in \mathbb{R}^q \ | \ \exists p_i\in P \ \text{s.t.} \  \|x-p_i\| \le \alpha \}.
\end{equation}
To ensure that we obtain a subcomplex of the Delaunay tessellation, we intersect $B_i(\alpha)$ with its corresponding Voronoi cell,
\begin{equation}
    R_i(\alpha) = B_i(\alpha) \cap V_i.
\end{equation}
It can be observed that $U(\alpha)=\cup_{p_i\in P}R_i(\alpha)$, so the $R_i$'s is a covering of $U(\alpha)$. The alpha complex $K_{\alpha}$ is the simplicial complex representing the nerve of this covering,
\begin{equation}
    K_{\alpha} = \{ J\subseteq \{1,2,...,|P|\} \;\; | \;\; \bigcap_{i\in J}  R_i(\alpha) \neq \emptyset \}.
\end{equation}
The equivalence to the original definition can be readily checked. The union of all simplices in the alpha complex forms the alpha shape.  \autoref{fig:Delaunay} illustrates the Voronoi diagram, Delaunay triangulation, and non-Delaunay triangulation. The point set is $P$ = \{A,B,C,D,E\}, and the blue lines in the left chart of \autoref{fig:Delaunay} separate the plane into  the Voronoi cells. The red circles are the empty circumcircles for triples of points in $P$. We can notice that no four points are on the same red circle, which satisfies the uniqueness condition for constructing the Delaunay triangulation. In the right chart of \autoref{fig:Delaunay}, the green circumcircle of ACD contains E and the green circumcirlce of AEC contains D, indicating that those two triangles do not belong to  the Delaunay triangulation.

\begin{figure}[ht!]
    \centering
		\includegraphics[width=0.8\textwidth]{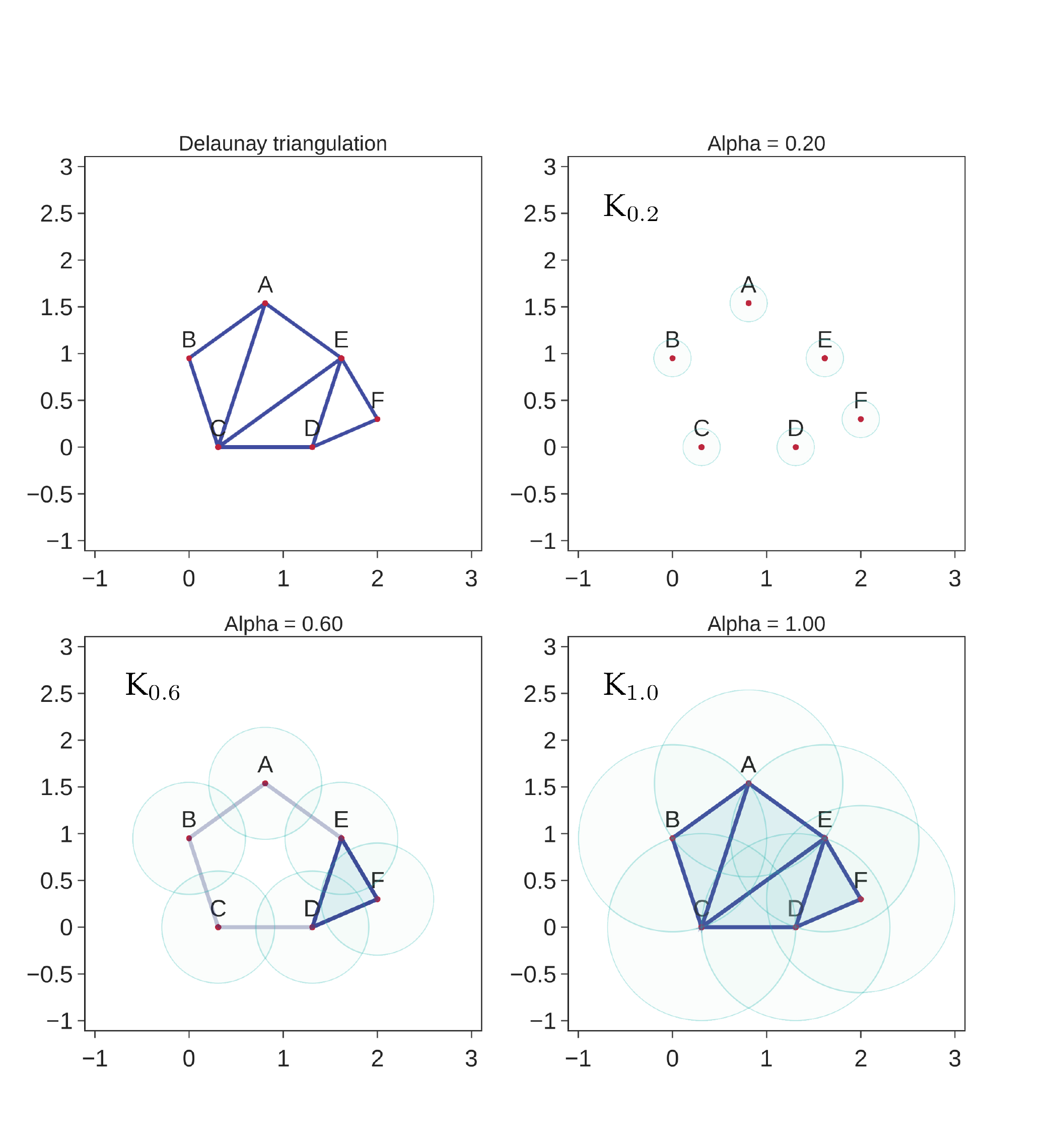}
    \caption{Illustration of 2D Delaunay triangulation, alpha shapes, and alpha complexes for a set of 6 points A, B, C, D, E, and F. {\bf Top left}: The 2D Delaunay triangulation. {\bf Top right}: The alpha shape and alpha complex at filtration value $\alpha = 0.2$. {\bf Bottom left }: The alpha shape and alpha complex at filtration value $\alpha = 0.6$. {\bf Bottom right}: The alpha shape and alpha complex at filtration value $\alpha = 1.0$. Here, we use dark blue color to fill the alpha shape.}
    \label{fig:Alpha2D}
\end{figure}

\begin{figure}[ht!]
    \centering
		\includegraphics[width=0.85\textwidth]{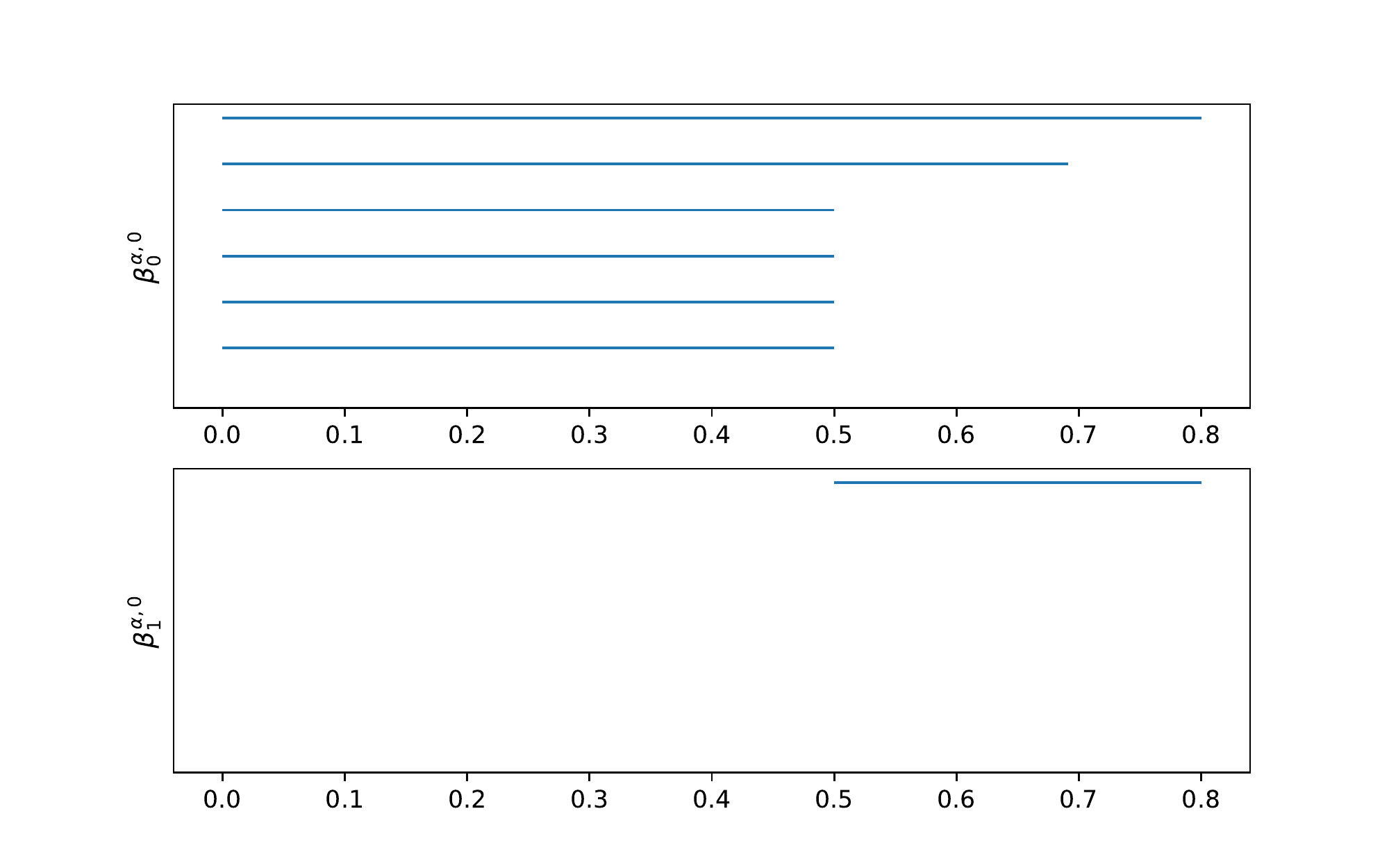}
    \caption{The persistent barcode for a set of points as illustrated in \autoref{fig:Alpha2D} that are generated from Gudhi and DioDe.}
    \label{fig:Betti2D}
\end{figure}

\autoref{fig:Alpha2D} illustrates the standard filtration of alpha complexes. The top left figure is the Delaunay triangulation of six 2D points A, B, C, D, E, and F. With an ever-growing radius $\alpha$ centered at these points, a family of sub-complexes of the Delaunay triangulation can be constructed. \autoref{fig:Betti2D} shows the persistence barcode of these 6 points. It can be seen that when $\alpha=0.2$, all six points are disconnected, indicating that 6 $0$-cycles (connected components) existed, which matches with \autoref{fig:Betti2D}, where there are a total of 6 bars when $\alpha=0.2$. With the radius $\alpha$ continually increasing, a $1$-cycle will be formed, and the associated alpha shape are shown in the bottom left chart of \autoref{fig:Alpha2D}. One can notice that in \autoref{fig:Betti2D}, when $\alpha = 0.6$, $\beta_{1}^{\alpha,0} = 1$. When $\alpha$ reaches $0.83$, the $1$-cycle disappears and $\beta_{1}^{\alpha,0} = 0$ as shown in the bottom left panel of \autoref{fig:Alpha2D}. \autoref{table:2DMatrix0} and \autoref{table:2DMatrix} show how we construct the $q$th-order persistent Laplacian $\mathcal{L}_{q}^{t,p}$ and calculate the harmonic ($\beta_{q}^{t,p}$) and non-harmonic persistent spectra of $\mathcal{L}_{q}^{t,p}$ from the simplicial complexes $K_{0.2}$ to $K_{0.6}$ and $K_{0.6}$ to $K_{0.6}$.

\begin{table}[th]
    \centering
    \scriptsize
    \setlength\tabcolsep{1pt}
    \captionsetup{margin=1cm}
    \caption{The matrix representation of $q$-boundary operator and its $q$th-order persistent  Laplacian with  corresponding dimension, rank, nullity, and spectra  from alpha complex $K_{0.6} \to K_{0.6}$.}
    \begin{tabular}{c|ccc}
    \hline
    $q$  & $q=0$ & $q=1$ & $q=2$ \\ \hline \\
    $\mathcal{B}_{q+1}^{0.6,0}$  & $\begin{array}{@{}r@{}c@{}c@{}c@{}c@{}c@{}c@{}c@{}c@{}l@{}}
            & \text{AB} & \text{BC} & \text{CD} & \text{DE} & \text{EF} & \text{DF} & \text{AE} \\
           \left.\begin{array}{c}
            \text{A} \\
            \text{B} \\
            \text{C} \\
            \text{D} \\
            \text{E} \\
            \text{F}
            \end{array}\right[
            & \begin{array}{c} -1 \\  1  \\  0  \\  0 \\ 0 \\ 0   \end{array}
            & \begin{array}{c} 0 \\   -1  \\  1  \\  0 \\ 0 \\0   \end{array}
            & \begin{array}{c} 0 \\   0  \\  -1  \\  1 \\ 0 \\0   \end{array}
            & \begin{array}{c} 0 \\  0  \\  0  \\  -1 \\ 1 \\0   \end{array}
            & \begin{array}{c} 0 \\  0  \\  0  \\  0 \\ -1 \\1   \end{array}
            & \begin{array}{c} 0 \\  0  \\  0  \\  -1 \\ 0 \\1  \end{array}
            & \begin{array}{c} -1 \\  0  \\  0  \\  0 \\ 1 \\0   \end{array}
            & \left]\begin{array}{c} \\ \\  \\ \\ \\ \\
            \end{array}
            \right.\end{array}$ & $\begin{array}{@{}r@{}c@{}c@{}c@{}c@{}c@{}c@{}c@{}c@{}l@{}}
                & \text{DEF} \\
           \left.\begin{array}{c}
            \text{AB} \\
            \text{BC} \\
            \text{CD} \\
            \text{DE} \\
            \text{EF} \\
            \text{DF} \\
            \text{AE}
            \end{array}\right[
            & \begin{array}{c} 0 \\  0  \\  0  \\  1 \\ 1 \\-1 \\ 0 \end{array}
            & \left]\begin{array}{c} \\ \\  \\ \\ \\ \\ \\
            \end{array}
            \right.\end{array}$ &  /    \\  \\
    $\mathcal{B}_{q}^{0.6}$  &  $\begin{array}{@{}r@{}c@{}c@{}c@{}c@{}c@{}c@{}l@{}}
            & \text{A} & \text{B} & \text{C} & \text{D} & \text{E} & \text{F}   \\
           \left.\begin{array}{c}
            \end{array}\right[
            & \begin{array}{c}  0  \end{array}
            & \begin{array}{c}  0  \end{array}
            & \begin{array}{c}  0  \end{array}
            & \begin{array}{c}  0  \end{array}
            & \begin{array}{c}  0  \end{array}
            & \begin{array}{c}  0  \end{array}
            & \left]\begin{array}{c}
            \end{array}\right.
        \end{array}$ & $\begin{array}{@{}r@{}c@{}c@{}c@{}c@{}c@{}c@{}c@{}c@{}l@{}}
            & \text{AB} & \text{BC} & \text{CD} & \text{DE} & \text{EF} & \text{DF} & \text{AE} \\
           \left.\begin{array}{c}
            \text{A} \\
            \text{B} \\
            \text{C} \\
            \text{D} \\
            \text{E} \\
            \text{F}
            \end{array}\right[
            & \begin{array}{c} -1 \\  1  \\  0  \\  0 \\ 0 \\ 0   \end{array}
            & \begin{array}{c} 0 \\   -1  \\  1  \\  0 \\ 0 \\0   \end{array}
            & \begin{array}{c} 0 \\   0  \\  -1  \\  1 \\ 0 \\0   \end{array}
            & \begin{array}{c} 0 \\  0  \\  0  \\  -1 \\ 1 \\0   \end{array}
            & \begin{array}{c} 0 \\  0  \\  0  \\  0 \\ -1 \\1   \end{array}
            & \begin{array}{c} 0 \\  0  \\  0  \\  -1 \\ 0 \\1  \end{array}
            & \begin{array}{c} -1 \\  0  \\  0  \\  0 \\ 1 \\0   \end{array}
            & \left]\begin{array}{c} \\ \\  \\ \\ \\ \\
            \end{array}
            \right.\end{array}$ &  $\begin{array}{@{}r@{}c@{}c@{}c@{}c@{}c@{}c@{}l@{}}
                & \text{DEF} \\
           \left.\begin{array}{c}
            \text{AB} \\
            \text{BC} \\
            \text{CD} \\
            \text{DE} \\
            \text{EF} \\
            \text{DF} \\
            \text{AE}
            \end{array}\right[
            & \begin{array}{c} 0 \\  0  \\  0  \\  1 \\ 1 \\-1\\0 \end{array}
            & \left]\begin{array}{c} \\ \\  \\ \\ \\ \\ \\
            \end{array}
            \right.\end{array}$    \\  \\
    $\mathcal{L}_{q}^{0.6,0}$  & $\left[\begin{array}{ccccccc}
         2  &  -1  &   0  &   0  &  -1 &  0  \\
        -1  &   2  &  -1  &   0  &  0 &  0 \\
         0  &  -1  &   2  &  -1  &  0 &  0  \\
         0  &   0  &  -1  &   3  & -1 & -1 \\
         -1  &   0  &   0  &  -1  &  3 & -1 \\
         0  &   0  &   0  &  -1  & -1 &  2
    \end{array}\right]$  & $\left[\begin{array}{ccccccccc}
         2  &  -1  &   0  &   0  &  0 &  0 & 1  \\
        -1  &   2  &  -1  &   0  &  0 &  0 & 0  \\
         0  &  -1  &   2  &  -1  &  0 & -1 & 0  \\
         0  &   0  &  -1  &   3  &  0 &  0 & 1  \\
         0  &   0  &   0  &   0  &  3 &  0 & -1 \\
         0  &   0  &  -1  &   0  &  0 &  3 & 0  \\
         1  &   0  &   0  &   1  & -1 &  0 & 2
    \end{array}\right]$ &  [3]    \\  \\
    $\beta_{q}^{0.6,0}$                        & 1               & 1 & 0    \\  \\
    $\text{dim}(\mathcal{L}_{q}^{0.6,0})$      & 6               & 7 & 1    \\  \\
    $\text{rank}(\mathcal{L}_{q}^{0.6,0})$     & 5               & 6 & 1    \\  \\
    $\text{nullity}(\mathcal{L}_{q}^{0.6,0})$  & 1               & 1 & 0  \\  \\
    $\text{Spec}(\mathcal{L}_{q}^{0.6,0})$    & $\{0,1, 1.5858,3,4,4.4142\}$   & $\{0,1,1.5858,3,3,4,4.4142\}$ & $\{3\}$  \\ \hline
    \end{tabular}
    \\
    \label{table:2DMatrix0}
\end{table}

\begin{table}[th]
    \centering
    \scriptsize
    \setlength\tabcolsep{10pt}
    \captionsetup{margin=1cm}
    \caption{The matrix representation of $q$-boundary operator and its  $q$th-order persistent  Laplacian with  corresponding dimension, rank, nullity, and spectra  from alpha complex $K_{0.2} \to K_{0.6}$.}
    \begin{tabular}{c|ccc}
    \hline
    $q$  & $q=0$ & $q=1$ & $q=2$ \\ \hline \\
    $\mathcal{B}_{q+1}^{0.2,0.4}$  & $\begin{array}{@{}r@{}c@{}c@{}c@{}c@{}c@{}c@{}c@{}c@{}l@{}}
            & \text{AB} & \text{BC} & \text{CD} & \text{DE} & \text{EF} & \text{DF} & \text{AE} \\
           \left.\begin{array}{c}
            \text{A} \\
            \text{B} \\
            \text{C} \\
            \text{D} \\
            \text{E} \\
            \text{F}
            \end{array}\right[
            & \begin{array}{c} -1 \\  1  \\  0  \\  0 \\ 0 \\ 0   \end{array}
            & \begin{array}{c} 0 \\   -1  \\  1  \\  0 \\ 0 \\0   \end{array}
            & \begin{array}{c} 0 \\   0  \\  -1  \\  1 \\ 0 \\0   \end{array}
            & \begin{array}{c} 0 \\  0  \\  0  \\  -1 \\ 1 \\0   \end{array}
            & \begin{array}{c} 0 \\  0  \\  0  \\  0 \\ -1 \\1   \end{array}
            & \begin{array}{c} 0 \\  0  \\  0  \\  -1 \\ 0 \\1  \end{array}
            & \begin{array}{c} -1 \\  0  \\  0  \\  0 \\ 1 \\0   \end{array}
            & \left]\begin{array}{c} \\ \\  \\ \\ \\ \\
            \end{array}
            \right.\end{array}$ & / &  /    \\  \\
    $\mathcal{B}_{q}^{0.2}$  &  $\begin{array}{@{}r@{}c@{}c@{}c@{}c@{}c@{}c@{}l@{}}
            & \text{A} & \text{B} & \text{C} & \text{D} & \text{E} & \text{F}   \\
           \left.\begin{array}{c}
            \end{array}\right[
            & \begin{array}{c}  0  \end{array}
            & \begin{array}{c}  0  \end{array}
            & \begin{array}{c}  0  \end{array}
            & \begin{array}{c}  0  \end{array}
            & \begin{array}{c}  0  \end{array}
            & \begin{array}{c}  0  \end{array}
            & \left]\begin{array}{c}
            \end{array}\right.
        \end{array}$ & / &  /    \\  \\
    $\mathcal{L}_{q}^{0.2,0.4}$  & $\left[\begin{array}{ccccccc}
         2  &  -1  &   0  &   0  &  -1 &  0  \\
        -1  &   2  &  -1  &   0  &  0 &  0 \\
         0  &  -1  &   2  &  -1  &  0 &  0  \\
         0  &   0  &  -1  &   3  & -1 & -1 \\
         -1  &   0  &   0  &  -1  &  3 & -1 \\
         0  &   0  &   0  &  -1  & -1 &  2
    \end{array}\right]$  & / &  /    \\  \\
    $\beta_{q}^{0.2,0.4}$                        & 1               & / & /    \\  \\
    $\text{dim}(\mathcal{L}_{q}^{0.2,0.4})$      & 6               & / & /    \\  \\
    $\text{rank}(\mathcal{L}_{q}^{0.2,0.4})$     & 5               & / & /    \\  \\
    $\text{nullity}(\mathcal{L}_{q}^{0.2,0.4})$  & 1               & / & /  \\  \\
    $\text{Spec}(\mathcal{L}_{q}^{0.2,0.4})$    & $\{0,1, 1.5858,3,4,4.4142\}$   & / & /  \\ \hline
    \end{tabular}
    \\
    \label{table:2DMatrix}
\end{table}

\section{Implementation}

\subsection{Construction of alpha shape}

Recall that, given a set of points, the alpha shape with any $\alpha$ value is a subcomplex of Delaunay tessellation. Thus, to construct the filtration of alpha complexes, it is necessary to first compute the complete simplicial complex through the Delaunay tessellation formed by the set of points. A number of efficient implementations is available in existing software packages. Our implementation employs the Computational Geometry Algorithms Library (CGAL), an efficient and robust software package for many commonly used calculations. We then assign each simplex $\sigma$ with an alpha value $\alpha_\sigma$. Finally, the alpha shape given at an $\alpha$ value $\alpha_0$ is constructed by union of convex hulls of all the simplices $\sigma$ satisfying $\alpha_\sigma \leq \alpha_0$, which naturally forms the nerve of balls centered at the given points truncated by the Voronoi regions, i.e., the corresponding alpha complex.

We illustrate our implementation with point sets $P$ in 3D, as it is the most common use scenario. We also assume that all the points are in general positions, which means that no 4 points of $P$ lie on the same plane and no 5 points of $P$ lie on the same sphere. Given a simplex $\sigma$, which can be a point, an edge, a triangle or a tetrahedron, denote the open ball bounded by its minimal circumsphere as $B_\sigma$. The simplex $\sigma$ is called \textit{Gabriel} (\cite{kerber2012medusa}) if $ B_\sigma \cap P = \emptyset.$ Note that for vertices (0-simplices) the circumradius is considered $0$. The above discussion can be directly adapted for 2D implementation by replacing circumsphere with circumcircle and omitting tetrahedra.

The filtration parameter $\alpha$ for every simplex $\sigma$ can be defined as follows. If the simplex is Gabriel, the filtration value is the corresponding circumradius (for efficiency, we actually store its square) because the corresponding ball can be considered as an empty $\alpha$-ball touching all its vertices. If the simplex is not Gabriel, the filtration value is the minimum of all the filtration values of the cofaces of $\sigma$ that contain the points making the simplex non-Gabriel. When $\alpha$ value reaches that number, we will have an empty $\alpha$-ball making the simplex $\alpha$-exposed.

\subsection{Implementation details}
To ensure the valid calculation of the filtration parameter for non-Gabriel simplices, the filtration value are always computed from the highest dimension (tetrahedra) down to 0 (vertices). We initialize the filtration value for all the simplices to be positive infinity. For dimension $k$, we iterate through each $k$-simplex. If the current filtration value $\alpha_\sigma^2$ is positive infinity, we assign the filtration value as the square of the corresponding circumradius. Then, we check every $(k\!-\!1)$-dimensional face $\tau$ in $\partial \sigma$. If the circumsphere of $\tau$ enclosed the other vertex of $\sigma$ in the interior, it is not Gabriel, and does not correspond to an empty $\alpha$-ball. In this case, $\alpha_\sigma^2$ is assigned to $\alpha_\tau^2$ if $\alpha_\sigma > \alpha_\tau.$

With this procedure, we ensure that $\alpha_\sigma$ for every simplex $\sigma$ corresponding to the filtration value $\alpha$ is $\alpha$-exposed to an empty $\alpha$-ball. In other words, we ensure that for each simplex represented by its vertex index set $J\subseteq \{1,2,...,|P|\}$ is in the nerve of $R_i$'s, which are the intersections $R_i=V_i\cap B_i$ of Voronoi cells $V_i$'s and balls $B_i$'s around the points $p_i$'s.

\subsubsection{Boundary operator construction}
With $\alpha_\sigma$ assigned, we sort the $k$-simplices with increasing filtration parameter value. This allows us to construct a single boundary operator $B_q^\infty$ (the matrix representation of $\partial_q^\infty$) for the entire filtration, which is that of the Delaunay tessellation. For any given $\alpha$, we can read of the top left block of the full boundary matrix $B_q^\infty$, i.e.,
\begin{equation}
    \left(B_q^\alpha\right)_{ij}=\left(B_q^\infty\right)_{ij},\quad \forall 1\leq i \leq N_{q-1}^\alpha, 1\leq j \leq N_{q}^\alpha,
\end{equation}
where $N_q^\alpha$ is the number of $q$-simplices in the alpha complex with the filtration parameter $\alpha.$ Alternative, we can consider the $N_q^\alpha\!\times\! N_q^\infty$ projection matrix $P_q^\alpha$ from the Delaunay tessellation to the alpha complex, $\left(P_q^\alpha\right)_{ij}=\delta_{ij}$ ($1$ on the diagonal and $0$ elsewhere), with which we have $B_q^\alpha = P_{q-1}^\alpha B_q^\infty (P_q^\alpha)^T.$

\subsubsection{Persistent boundary operator} The construction of $p$-persistent boundary matrix $B_q^{\alpha,p}$ (the representation of operator $\eth_q^{\alpha,p}$ is more involved than reading off $B_q^\infty.$

We first construct the projection matrix $\mathbb{P}_q^{\alpha,p}$ from $C_q^{\alpha+p}$ to $\mathbb{C}_q^{\alpha,p}.$ Then, the $p$-persistent boundary matrix can be assembled as $B_q^{\alpha,p}= P_{q-1}^\alpha  B_q^\infty (\mathbb{P}_{q}^{\alpha,p})^T.$

To construct the projection matrix, we first note that it is the projection to the kernel of an operator that measures the difference between the boundary operator mapped onto $C_{q-1}^{\alpha+p}$ and the boundary restricted to $C_{q-1}^\alpha$, $\text{Diff}_{q}^{\alpha,p}=(I_{q-1}^{\alpha+p}-R_{q-1}^{\alpha,p})^T B_{q}^{\alpha+p},$ where $R_q^{\alpha,p}=P_q^{\alpha+p} (P_q^\alpha)^T P_q^\alpha (P_q^{\alpha+p})^T$ is the restriction from $C_q^{\alpha+p}$ to $C_q^\alpha$ and $I_q^{\alpha+p}$ is the identity matrix on $C_q^{\alpha+p}$.

Instead of storing a dense matrix, we propose to use a procedural representation involving the inverse of persistent Laplacians with gauge (\cite{zhao20193d}) to reduce the storage as well as speed up the computation. More specifically, we construct the projection matrix as follows
\begin{equation}
\mathbb{P}_{q}^{\alpha,p}=I_{q}^{\alpha+p}-(\tilde{\text{Diff}}_{q}^{\alpha,p})^T(\tilde{L}_{q-1}^{\alpha,p})^{-1} \tilde{\text{Diff}}_{q}^{\alpha,p},
\end{equation}
where $(\tilde{L}_{q-1}^{\alpha,p})^{-1}$ can be implemented through rank deficiency fixing in \cite{zhao20193d}, and the restricted operator $\tilde{\text{Diff}}_{q}^{\alpha,p}$ is defined below. Note that this sparse linear equation solving approach is essentially the graph version of the harmonic extension described in Ref. \cite{zhao2020rham}.

The reason that the projection matrix can be defined this way is that starting from an arbitrary element $\omega_q \in C_q^{\alpha+p}$, we can modify it into $\omega_q - (\text{Diff}_{q}^{\alpha,p})^T f_{q-1} \in \mathbb{C}_q^{\alpha,p},$ where $f_{q-1}$ is nonzero only in the difference complex $\text{Cl}(T_{\alpha+p}-T_\alpha)$, the closure of the difference between $T_{\alpha+p}$ and $T_{\alpha}$. Denoting any chain $f$ on the difference complex as $\tilde{f}$ and any operator $B$ on it as $\tilde{B}^{\alpha,p},$ and the  $\tilde{B}_{q}^{\alpha,p} (\tilde{B}_{q}^{\alpha,p})^T \tilde{f}_{q-1} = \tilde{B}_{q}^{\alpha,p} \tilde{\omega}_q.$ Noticing that $\tilde{f}_{q-1}$ is determined up to a gauge transform $f_{q-1}-(\tilde{B}_{q-1}^{\alpha,p})^T \tilde{g}_{q-2}$ for some $(q-2)$-chain $g_{q-2}$ in $\text{Cl}(T_{\alpha+p}-T_\alpha),$ we introduce the gauge fixing term $\tilde{B}_{q-1}^{\alpha,p} f_{q-1}=0$, which leads us to the sparse linear system $\tilde{L}_{q-1}^{\alpha+p} \tilde{f}_{q-1} =  \tilde{\text{Diff}}_{q}^{\alpha,p}\omega_q$ where the $\tilde{\text{Diff}}$ operator is the above operator projected to the difference complex. Note that fixing the rank deficiency of persistent Laplacians (in the difference complex) is computationally efficient as its kernel dimension is far smaller than that of the corresponding boundary or coboundary operators.

\subsubsection{Persistent spectrum computation}
The $q$-order $p$-persistent Laplacian operators can then be implemented by direct evaluation of $L_q^{\alpha,p}=B_{q+1}^{\alpha,p} (B_{q+1}^{\alpha,p} )^T+(B_q^\alpha)^T B_q^\alpha.$ Their spectra can be evaluated through any off-the-shelf sparse matrix eigensolver.

Thus, the dimension of the null space of $L_0^{\alpha,p}$ is number of $p$-persistent connected components. The dimension of the null space of $L_1^{\alpha,p}$ is number of $p$-persistent handles or tunnels. Similarly, the dimension of the null space of $L_2^{\alpha,p}$ is the number of $p$-persistent cavities.

\section{Validation}

We construct the alpha complex at different filtration values from the finite cells of a Delaunay tessellation from the Computational Geometry Algorithms Library (\href{https://doc.cgal.org/latest/Manual/packages.html}{CGAL}). Gudhi and DioDe are two of the most frequently applied open-source libraries that are able to compute the Betti numbers (harmonic persistent spectra) based on CGAL. As shown in \cite{wang2020persistent}, the $0$-persistent $q$th Betti numbers $\beta_{q}^{\alpha,0}$ at filtration parameter $\alpha$ is the number of zero eigenvalues of $q$th-order $0$-persistent Laplacian $\mathcal{L}_{q}^{\alpha,0}$:
\begin{equation}
    \beta_{q}^{\alpha,0} = \dim(C_q^{\alpha}) - \text{rank}(\mathcal{L}_{q}^{\alpha,0}) = \dim \ker \mathcal{L}_{q}^{\alpha,0}.
\end{equation}
In fact, $\beta_{q}^{\alpha,0}$ counts the number of $q$-cycles in alpha complex $K_{\alpha}$ that persists in $K_{\alpha}$. Although Gudhi and DioDe can calculate the number of zero eigenvalues, the non-harmonic persistent spectra also play an important role in applications as shown in our earlier work \cite{wang2020persistent}. Therefore, we developed an open-source package \emph{HERMES}, which not only tracks the topological changes from the persistent Betti numbers but also derives the geometric changes from the non-harmonic  spectra of  persistent Laplacians. In the following, we compare the Betti numbers $\beta_{q}^{t,p}$ that are calculated from HERMES with the Betti numbers that are derived from Gudhi and DioDe on a set of 2D and 3D points, aiming to validate the robustness and accuracy of HERMES.


\clearpage
\begin{figure}[ht!]
    \centering
		\includegraphics[width=1\textwidth]{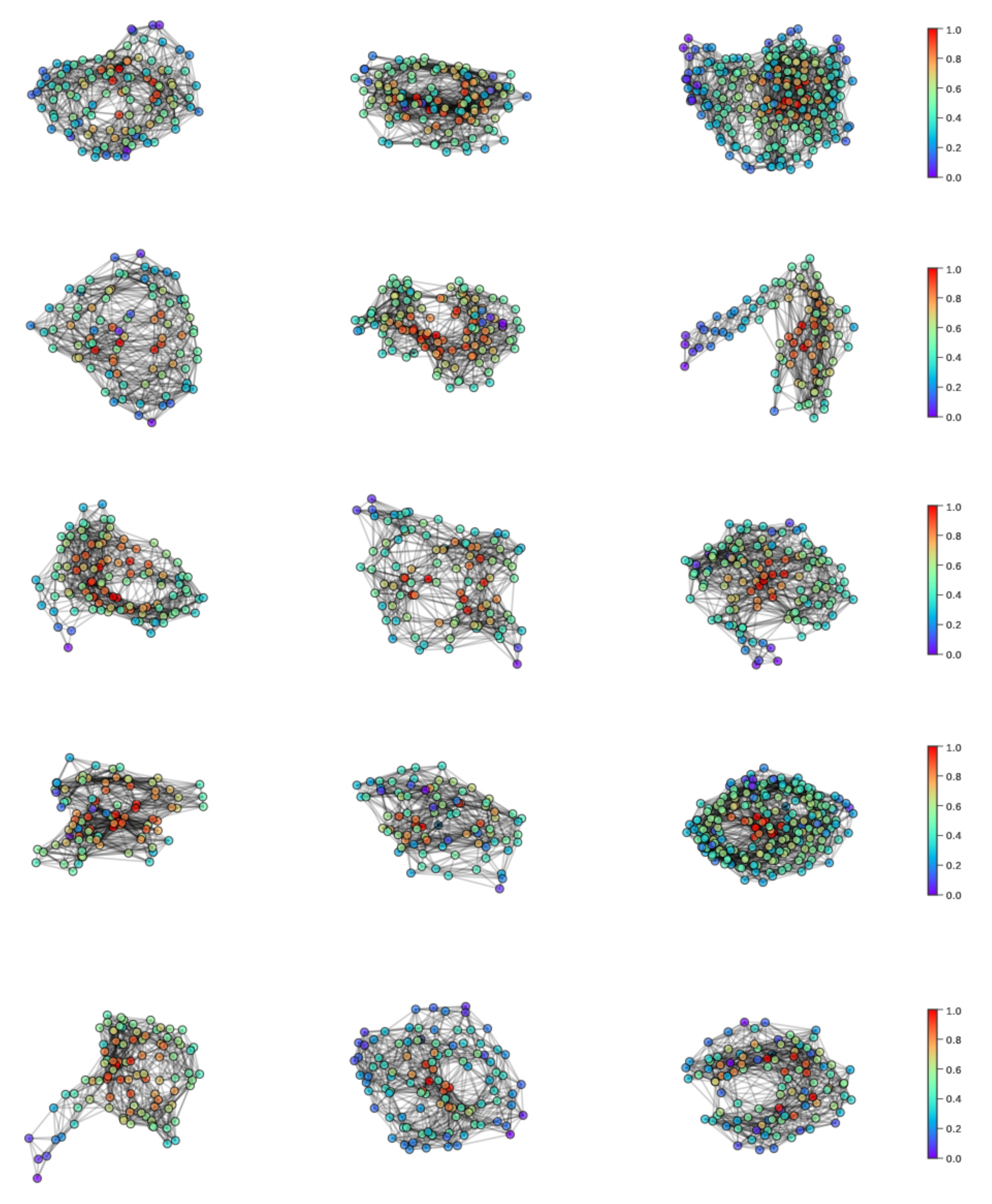}
    \caption{The alpha carbon network plots of 15 proteins: PDB IDs 1CCR, 1NKO, 1O08, 1OPD, 1QTO, 1R7J, 1V70, 1W2L, 1WHI, 2CG7, 2FQ3, 2HQK, 2PKT, 2VIM, and 5CYT from left to right and top to bottom. The color represents  the normalized diagonal element of the accumulated Laplacian  at each alpha carbon atom.}
    \label{fig:Network}
\end{figure}

\begin{figure}[ht!]
    \centering
		\includegraphics[width=0.8\textwidth]{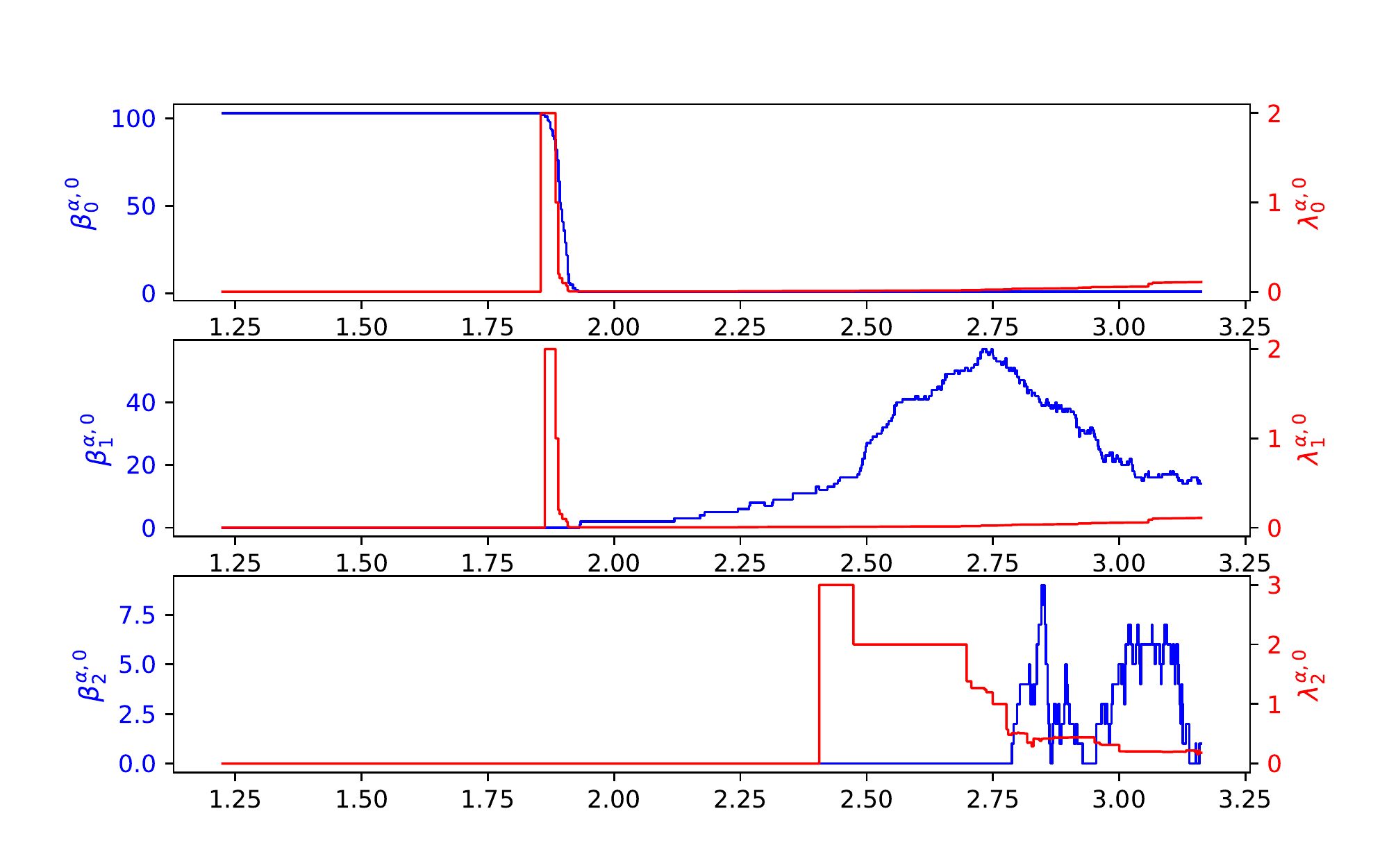}
    \caption{Illustration of the  harmonic spectra $\beta_q^{\alpha,0}$   (blue curve) and the smallest non-zero eigenvalue $\lambda_q^{\alpha,0}$ (red curve) of PDB ID 5CYT (the bottom left chart in Fig. \ref{fig:Network}) at different filtration value $\alpha$ when $q = 0,1,2$. The $\beta_q^{\alpha,0}$ are calculated from Gudhi, DioDe, and HERMES, and $\lambda_q^{\alpha,0}$ are obtained only from HERMES. Here, the $x$-axis represents the radius filtration value $\alpha$ (unit: $\SI{}{\angstrom}$), the left-$y$-axis represents  the number of zero eigenvalues of $\mathcal{L}_q^{\alpha,0}$, and the right-$y$-axis represents  the first non-zero eigenvalue of $\mathcal{L}_q^{\alpha,0}$. Note that the harmonic spectra from three methods are indistinguishable.
	}
    \label{fig:ProStepfuc}
\end{figure}

We first confirmed the correctness of HERMES with simple systems such as an icosahedron and C$_{70}$ molecules with known persistent Betti numbers \cite{xia2014persistent}. Results of these samples computed from Gudhi, DioDe, and HERMES are identical and are omitted because they are the same as those in the literature  \cite{xia2014persistent}.
Then, we apply Gudhi, DioDe, and HERMES to calculate the persistent Betti numbers of 15 proteins. Their Protein Data Bank (PDB) IDs are 1CCR, 1NKO, 1O08, 1OPD, 1QTO, 1R7J, 1V70, 1W2L, 1WHI, 2CG7, 2FQ3, 2HQK, 2PKT, 2VIM, and 5CYT. The 3D structures of these 15 proteins can be downloaded from the PDB (\url{https://www.rcsb.org/}). Here, only the alpha carbon atoms are considered in our calculations. \autoref{fig:Network} illustrates the network structures of  15 proteins. For each protein, color at atomic positions represents the normalized diagonal values of the accumulated 0th-order 0-persistent Laplacians:  $\frac{1}{\max \left(\mathcal{L}_0^{0}\right)_{ii}}\left(\mathcal{L}_0^{0}\right)_{ii}$, where $\mathcal{L}_0^{0}=\sum_\alpha \mathcal{L}_0^{\alpha,0}$. Here, the filtration $\alpha$ is from $\sqrt{1.5}$ \AA\ to $\sqrt{10}$ \AA\ with the $\SI{0.01}{\angstrom}$ step size. \autoref{fig:ProStepfuc} depicts the persistent Betti numbers $\beta_q^{\alpha,0}$  (blue curve) of PDB ID 5CYT that are calculated from Gudhi, DioDe, and HERMES,  together with the smallest non-zero eigenvalue $\lambda_q^{\alpha,0}$ (red curve) that are obtained only from HERMES. 

It can be seen that all of these three packages return exactly the same persistent Betti numbers, suggesting that the calculation of our package HERMES is reliable.  Additionally, the values of smallest non-zero eigenvalues $\lambda_0^{\alpha,0}$ and $\lambda_1^{\alpha,0}$ increase around $\SI{1.86}{\angstrom}$, indicating the dramatic  topological changes  at this point. Similarly, with the increment of the $\alpha$, the curve of $\lambda_2^{\alpha,0}$  also records the topological and geometric changes at a specific filtration value. The use of non-harmonic spectra for biophysical modeling was described  in our earlier work 
\cite{wang2020persistent}. 

\begin{figure}[ht!]
    \centering
		\includegraphics[width=0.8\textwidth]{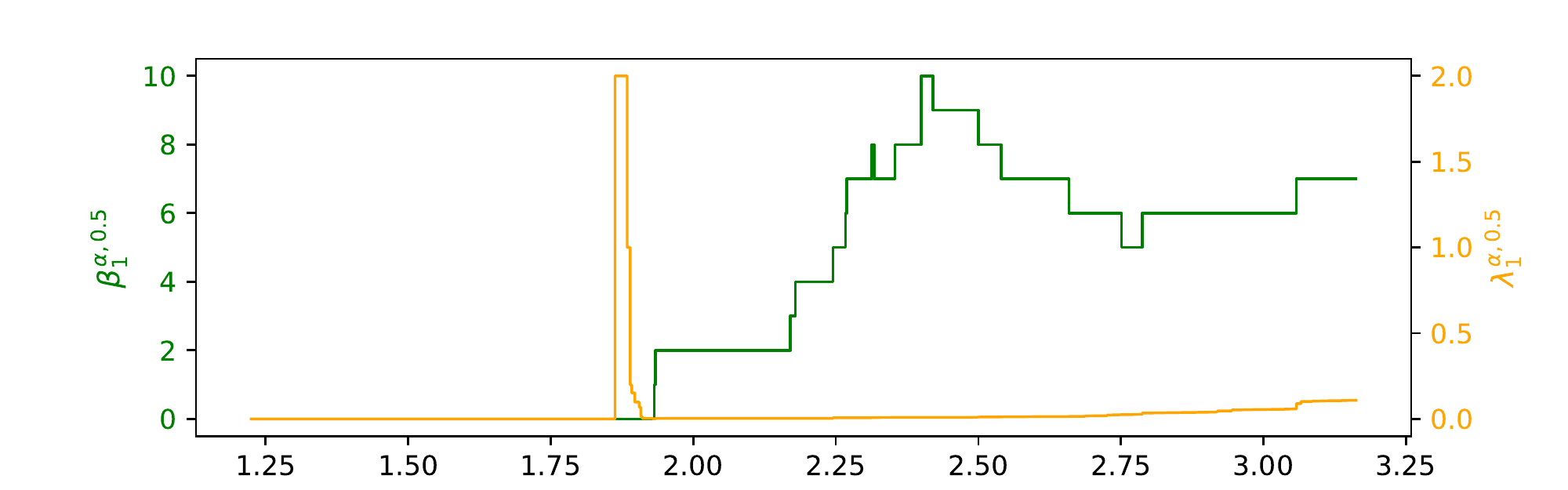}
    \caption{Illustration of the  harmonic spectra $\beta_1^{\alpha,0.5}$   (green curve) and the smallest non-zero eigenvalue $\lambda_1^{\alpha,0.5}$ (yellow curve) of PDB ID 5CYT (the bottom left chart in Fig. \ref{fig:Network}) at different filtration value $\alpha$ calculated from HERMES. Here, the $x$-axis represents the radius filtration value $\alpha$ (unit: $\SI{}{\angstrom}$), the left-$y$-axis represents the number of zero eigenvalues of $\mathcal{L}_1^{\alpha,0.5}$, and the right-$y$-axis represents  the first non-zero eigenvalue of $\mathcal{L}_1^{\alpha,0.5}$.
	}
    \label{fig:ProStepfuc_persistent}
\end{figure}

  To be noted, HERMES can also deal with the $q$th-order p-persistent Laplacians $\mathcal{L}_q^{\alpha,p}$.  \autoref{fig:ProStepfuc_persistent} illustrates the persistent Betti numbers $\beta_1^{\alpha,0.5}$ (green curve) and the smallest non-zero eigenvalue $\lambda_1^{\alpha,0.5}$ (yellow curve) of 5CYT that are computed from HERMES, demonstrating the capacity of HERMES for the direct calculation of the persistent spectra of $\mathcal{L}_q^{\alpha,p}$ ($p>0$). Compared with the middle chart of \autoref{fig:ProStepfuc}, the $\beta_1^{\alpha,0.5}$ in \autoref{fig:ProStepfuc_persistent} is always smaller than $\beta_1^{\alpha,0}$ at the same filtration $\alpha$. Moreover, the $\lambda_1^{\alpha,0.5}$ also goes up around $\SI{1.86}{\angstrom}$, which has the same behavior as $\lambda_1^{\alpha,0}$.

\begin{figure}[ht!]
    \centering
		\includegraphics[width=1\textwidth]{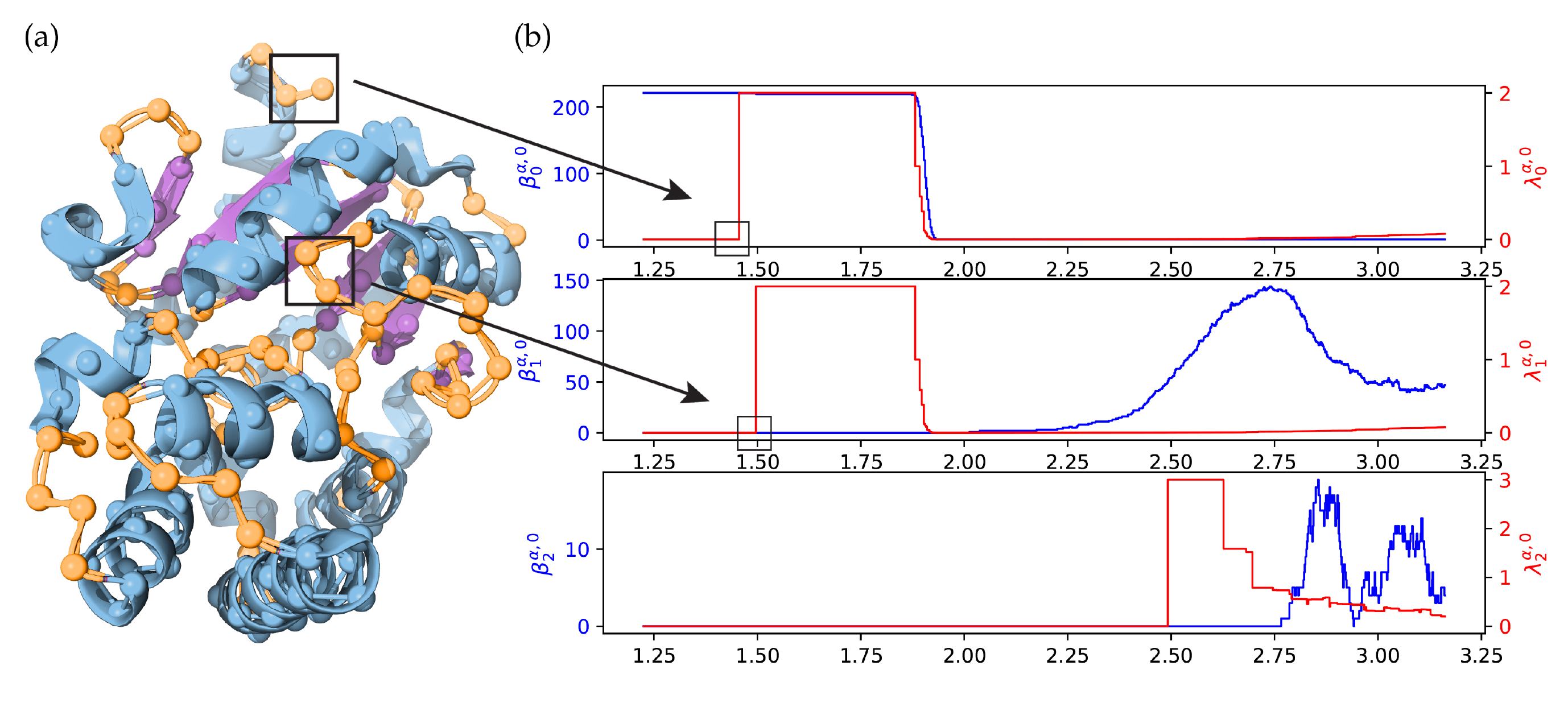}
    \caption{  (a) The 3D secondary structure of PDB ID 1O08. The blue, purple, and orange colors represent helix, sheet, and random coils of PDB ID 1O08. The ball represents the alpha carbon of PDB ID 1O08. (b) Illustration of the  harmonic spectra $\beta_q^{\alpha,0}$ (blue curve) and the smallest non-zero eigenvalue $\lambda_q^{\alpha,0}$ (red curve) of PDB ID 1O08  at different filtration value $\alpha$ when $q = 0,1,2$. The $\beta_q^{\alpha,0}$ are calculated from Gudhi, DioDe, and HERMES, and $\lambda_q^{\alpha,0}$ are calculated only from HERMES. Here, the $x$-axis represents the radius filtration value $\alpha$ (unit: $\SI{}{\angstrom}$), the left-$y$-axis represents for the number of zero eigenvalue of $\mathcal{L}_q^{\alpha,0}$, and the right-$y$-axis represents for the non-zero eigenvalues of $\mathcal{L}_q^{\alpha,0}$. Note that the harmonic spectra from three methods are indistinguishable. 
	}
    \label{fig:1O08_Combine}
\end{figure}

Furthermore, HERMES can be used to detect the abnormality of a protein structure. \autoref{fig:1O08_Combine} (a) shows a 3D secondary structure of PDB 1O08, where the balls represent the alpha carbon atoms. The light blue, purple, and orange colors represent helix, sheet, and random coils of PDB ID 1O08. \autoref{fig:1O08_Combine} (b) depicts its harmonic spectra $\beta_q^{\alpha,0}$   (blue curve) and the smallest non-zero eigenvalue $\lambda_q^{\alpha,0}$ (red curve). Notably,  two unusual onset of $\beta_0^{\alpha,0}$ and $\beta_1^{\alpha,0}$ are detected when $\alpha << \SI{1.9}{\angstrom}$, indicating something is wrong with the structure data. Usually, the distance between the two alpha carbon atoms is around $\SI{3.8}{\angstrom}$. By examining the structure of PDB 1O08, we found that two pairs of alpha carbon atoms in  PDB 1O08 have abnormal distances as  marked with black frames. The distance of alpha carbon atoms in the upper box is $\SI{2.914}{\angstrom}$ and that in the lower box is $\SI{2.996}{\angstrom}$, which are too short. 
 The plots of the other proteins can be found in the Appendix. Similar structural defects are detected for PDB IDs 1V70, 2HQK, 2PKT, and 2VIM,

Although our package provides additional geometric information by calculating the non-harmonic spectra of $q$th-order persistent Laplacians, there are two limitations of HERMES. First, HERMES is currently implemented only for the system of the alpha complexes. Additionally, the input format of HERMES is point cloud data. Other input formats, such as  pairwise distances, point cloud with van der Waals radii, and volumetric density are not supported. These limitations will be addressed in our future implementation.

\section{Conclusion}

While spectral graph theory has had tremendous success in data science to capture the geometric and topological information, it is limited by representing a graph structure at given  characteristic length scale, which hinders its practical application in data analysis. Motivated by the persistent (co)homology in dealing with the given initial data by constructing a family of simplicial complexes to track their topological invariants, and the multiscale graphs by creating a set of spectral graphs aiming to extract rich geometric information, we proposed persistent spectral graph (PSG) theory as a unified multiscale paradigm for simultaneous  geometric and topological analysis \cite{wang2019persistent}. PSG theory has stimulated mathematical analysis and algorithm development \cite{memoli2020persistent}, as well as applications to drug discovery \cite{meng2020persistent}, and protein flexibility analysis \cite{wang2020persistent}.

To enable broad and convenient applications of the PSG method, we present an open-source software package called  highly efficient robust multidimensional evolutionary spectra (HERMES). For a given point-cloud dataset, HERMES creates a persistent Laplacian (PL) at various topological dimensions via a filtration. The spectrum of PLs includes harmonic parts and non-harmonic parts. It turns out that the harmonic parts span the kernel spaces of PLs and carry the full topological information of the dataset. As a result, HERMES delivers the same topological data analysis (TDA) as does persistent homology. The non-harmonic parts of PLs provide valuable geometric analysis of the shape of data at various topological dimensions. The smallest non-zero eigenvalues are found to be very sensitive to  data abnormality. In the present HERMES implementation, only the alpha complex is realized while other complexes, such as Vietoris–Rips complex, will be made available in the near future. HERMES has been extensively validated for its accuracy, robustness, and reliability by  standard test datasets and a large number of protein structures.

\section*{Acknowledgment}
This work was supported in part by NIH grant  GM126189, NSF Grants DMS-1721024,  DMS-1761320, and IIS-1900473,  Michigan Economic Development Corporation,  George Mason University award PD45722,  Bristol-Myers Squibb, and Pfizer.


\begin{appendices}
\section{Supplementary figures}
\begin{figure}[ht!]
    \centering
		\includegraphics[width=0.8\textwidth]{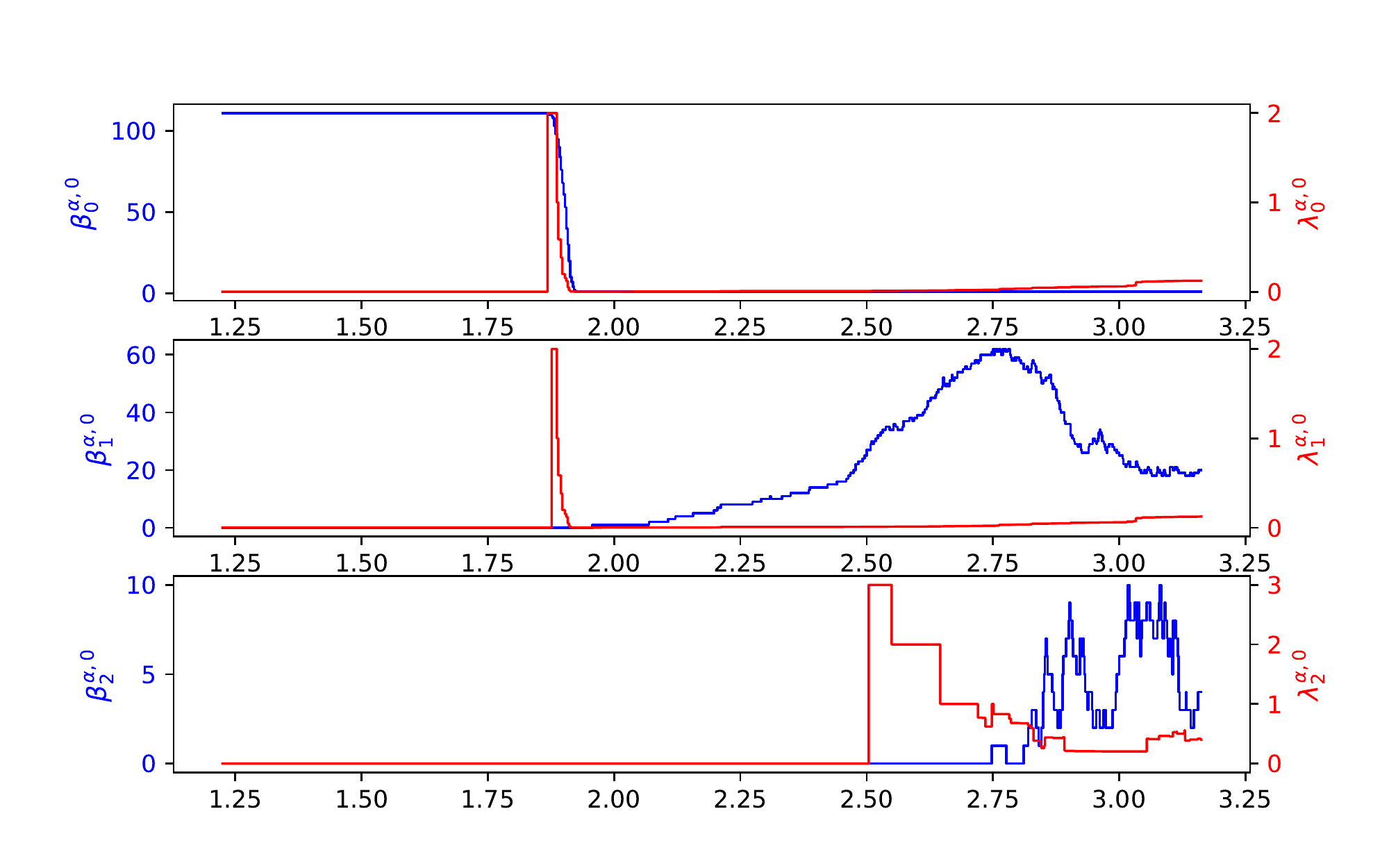}
    \caption{Illustration of the  harmonic spectra $\beta_q^{\alpha,0}$   (blue curve) and the smallest non-zero eigenvalue $\lambda_q^{\alpha,0}$ (red curve) of PDB ID 1CCR at different filtration value $\alpha$ when $q = 0,1,2$. The $\beta_q^{\alpha,0}$ are calculated from Gudhi, DioDe, and HERMES, and $\lambda_q^{\alpha,0}$ are obtained only from HERMES. Here, the $x$-axis represents the radius filtration value $\alpha$ (unit: $\SI{}{\angstrom}$), the left-$y$-axis represents   the number of zero eigenvalues of $\mathcal{L}_q^{\alpha,0}$, and the right-$y$-axis represents   the first non-zero eigenvalue of $\mathcal{L}_q^{\alpha,0}$. Note that the harmonic spectra from three methods are indistinguishable. 
		}
    \label{fig:1CCR}
\end{figure}

\begin{figure}[ht!]
    \centering
		\includegraphics[width=0.8\textwidth]{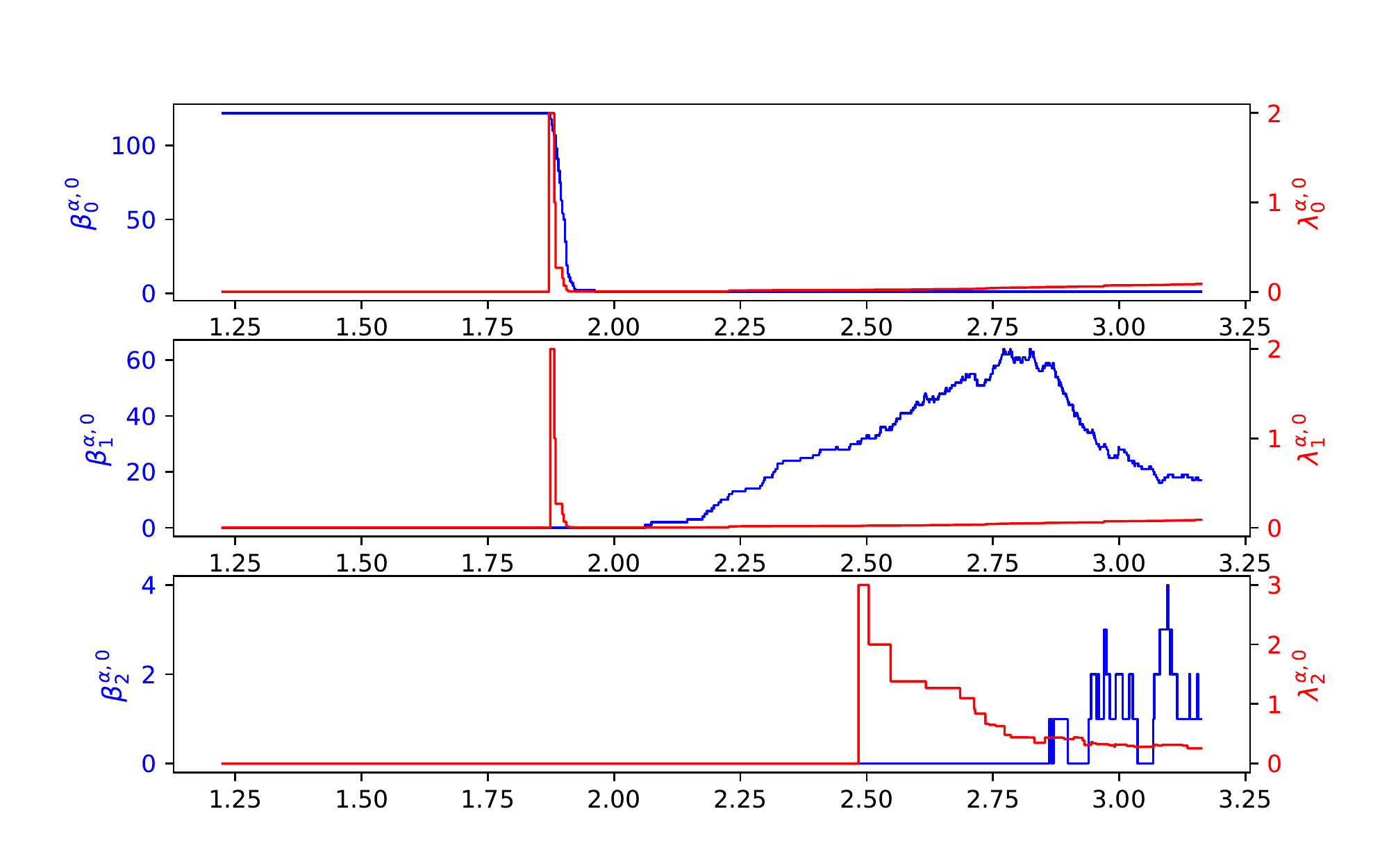}
    \caption{Illustration of the  harmonic spectra $\beta_q^{\alpha,0}$   (blue curve) and the smallest non-zero eigenvalue $\lambda_q^{\alpha,0}$ (red curve) of  PDB ID 1NKO at different filtration value $\alpha$ when $q = 0,1,2$. The $\beta_q^{\alpha,0}$ are calculated from Gudhi, DioDe, and HERMES, and $\lambda_q^{\alpha,0}$ are obtained only from HERMES. Here, the $x$-axis represents the radius filtration value $\alpha$ (unit: $\SI{}{\angstrom}$), the left-$y$-axis represents   the number of zero eigenvalues of $\mathcal{L}_q^{\alpha,0}$, and the right-$y$-axis represents   the first non-zero eigenvalue of $\mathcal{L}_q^{\alpha,0}$. Note that the harmonic spectra from three methods are indistinguishable. 
		}
    \label{fig:1NKO}
\end{figure}

\begin{figure}[ht!]
    \centering
		\includegraphics[width=0.8\textwidth]{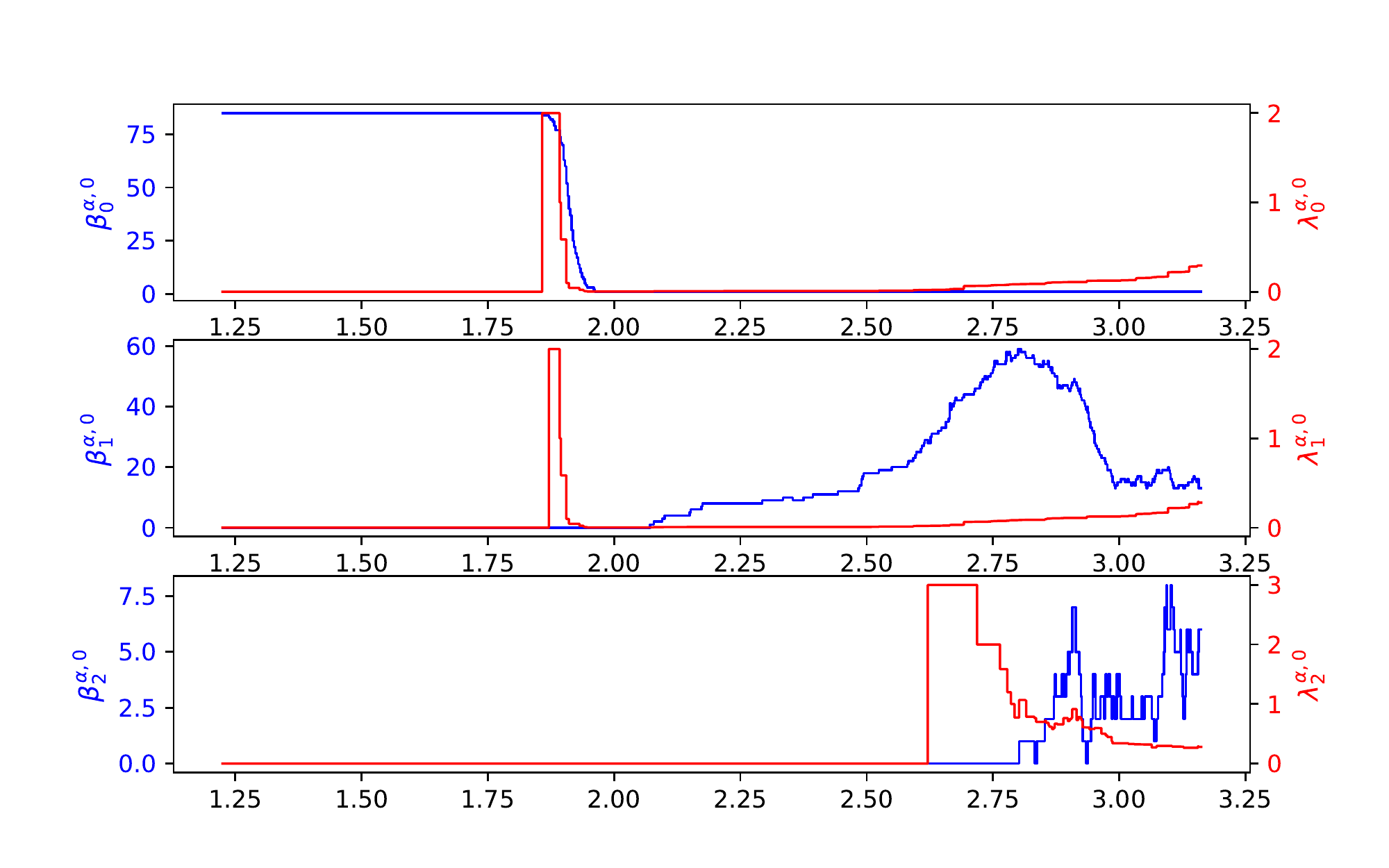}
    \caption{Illustration of the  harmonic spectra $\beta_q^{\alpha,0}$   (blue curve) and the smallest non-zero eigenvalue $\lambda_q^{\alpha,0}$ (red curve) of  PDB ID 1OPD at different filtration value $\alpha$ when $q = 0,1,2$. The $\beta_q^{\alpha,0}$ are calculated from Gudhi, DioDe, and HERMES, and $\lambda_q^{\alpha,0}$ are obtained only from HERMES. Here, the $x$-axis represents the radius filtration value $\alpha$ (unit: $\SI{}{\angstrom}$), the left-$y$-axis represents   the number of zero eigenvalues of $\mathcal{L}_q^{\alpha,0}$, and the right-$y$-axis represents   the first non-zero eigenvalue of $\mathcal{L}_q^{\alpha,0}$. Note that the harmonic spectra from three methods are indistinguishable. 
		}
    \label{fig:1OPD}
\end{figure}

\begin{figure}[ht!]
    \centering
		\includegraphics[width=0.8\textwidth]{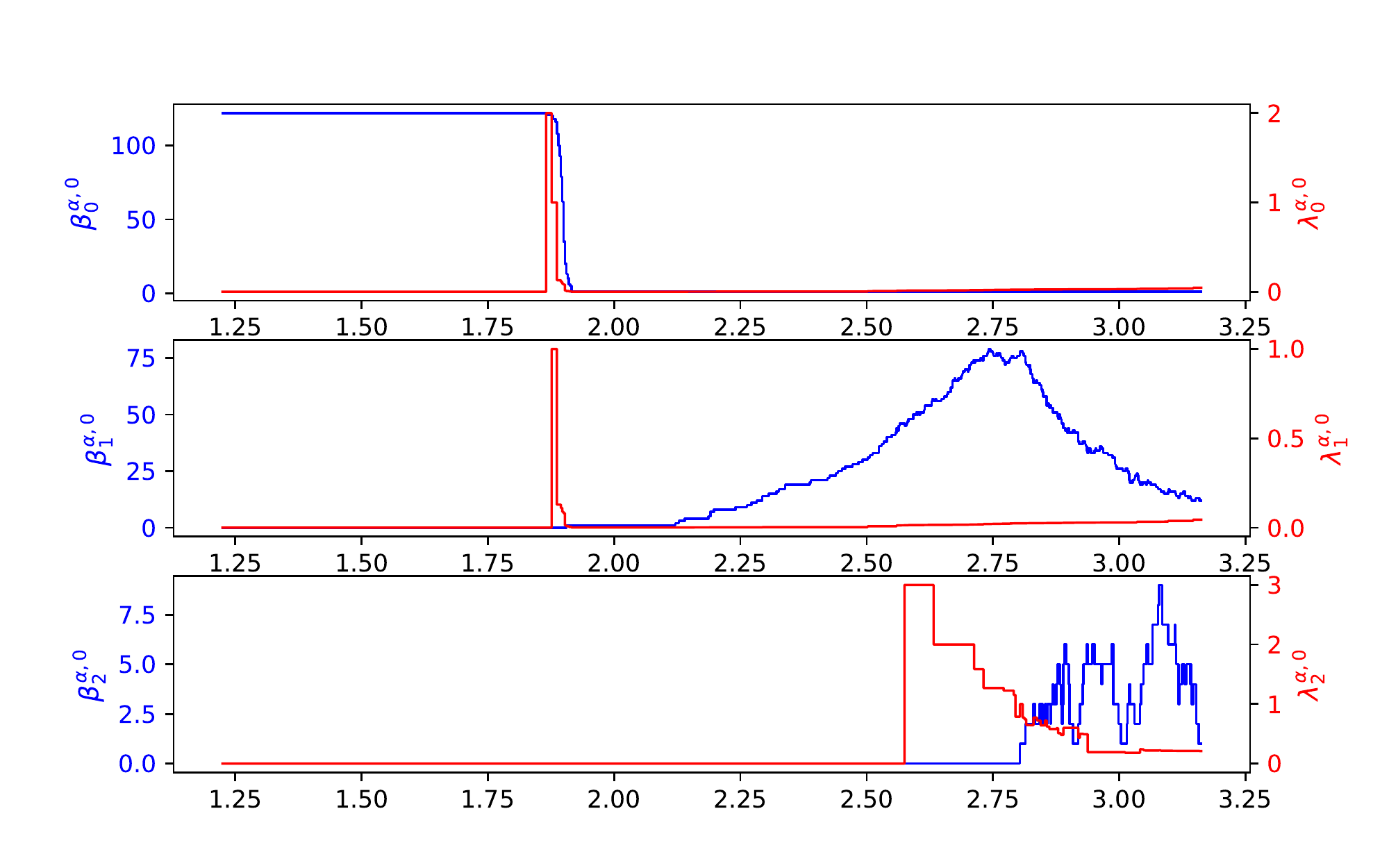}
    \caption{Illustration of the  harmonic spectra $\beta_q^{\alpha,0}$   (blue curve) and the smallest non-zero eigenvalue $\lambda_q^{\alpha,0}$ (red curve) of  PDB ID 1QTO at different filtration value $\alpha$ when $q = 0,1,2$. The $\beta_q^{\alpha,0}$ are calculated from Gudhi, DioDe, and HERMES, and $\lambda_q^{\alpha,0}$ are obtained only from HERMES. Here, the $x$-axis represents the radius filtration value $\alpha$ (unit: $\SI{}{\angstrom}$), the left-$y$-axis represents   the number of zero eigenvalues of $\mathcal{L}_q^{\alpha,0}$, and the right-$y$-axis represents   the first non-zero eigenvalue of $\mathcal{L}_q^{\alpha,0}$. Note that the harmonic spectra from three methods are indistinguishable. 
		}
    \label{fig:1QTO}
\end{figure}

\begin{figure}[ht!]
    \centering
		\includegraphics[width=0.8\textwidth]{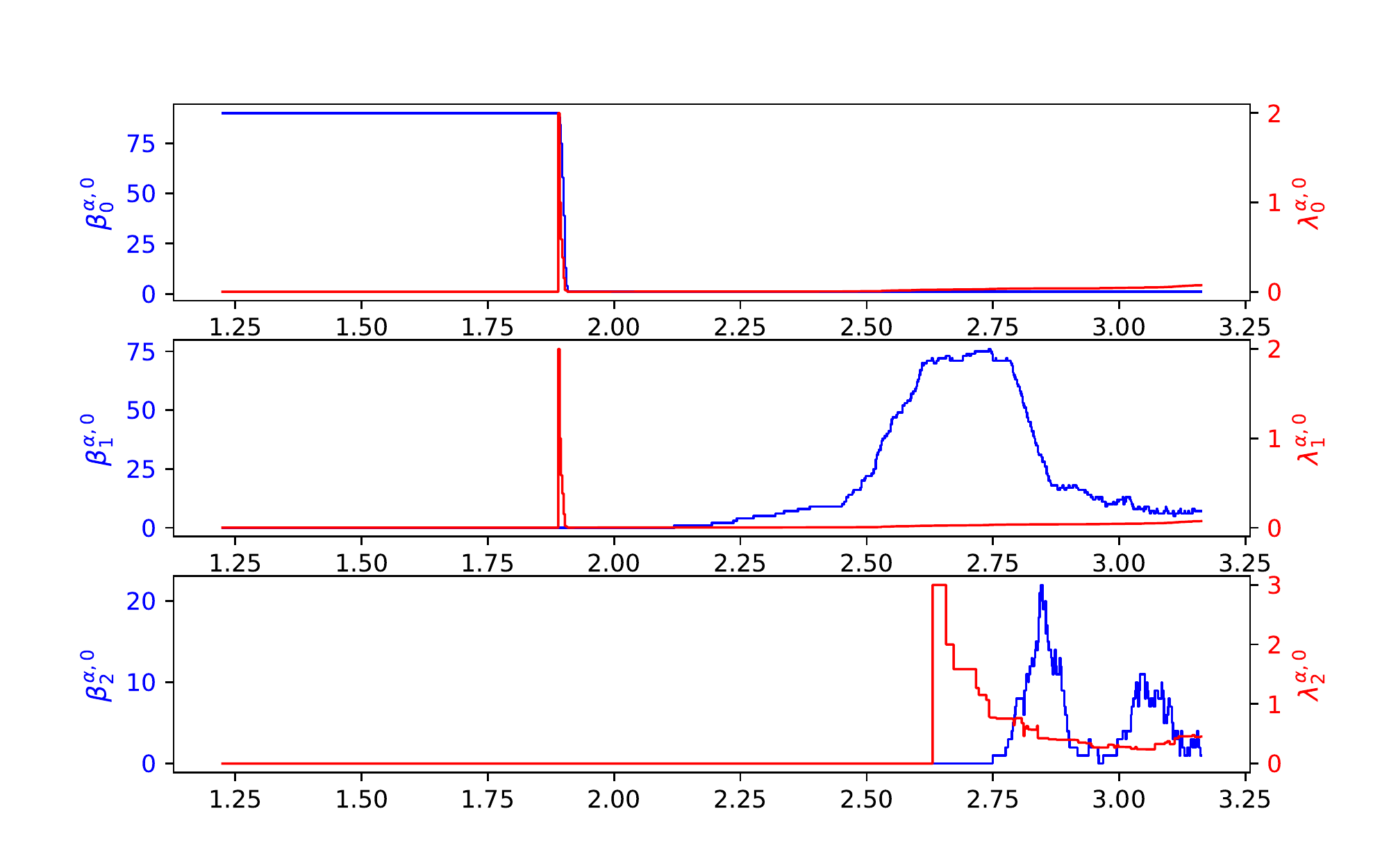}
    \caption{Illustration of the  harmonic spectra $\beta_q^{\alpha,0}$   (blue curve) and the smallest non-zero eigenvalue $\lambda_q^{\alpha,0}$ (red curve) of  PDB ID 1R7J at different filtration value $\alpha$ when $q = 0,1,2$. The $\beta_q^{\alpha,0}$ are calculated from Gudhi, DioDe, and HERMES, and $\lambda_q^{\alpha,0}$ are obtained only from HERMES. Here, the $x$-axis represents the radius filtration value $\alpha$ (unit: $\SI{}{\angstrom}$), the left-$y$-axis represents   the number of zero eigenvalues of $\mathcal{L}_q^{\alpha,0}$, and the right-$y$-axis represents   the first non-zero eigenvalue of $\mathcal{L}_q^{\alpha,0}$. Note that the harmonic spectra from three methods are indistinguishable.
		}
    \label{fig:1R7J}
\end{figure}

\begin{figure}[ht!]
    \centering
		\includegraphics[width=0.8\textwidth]{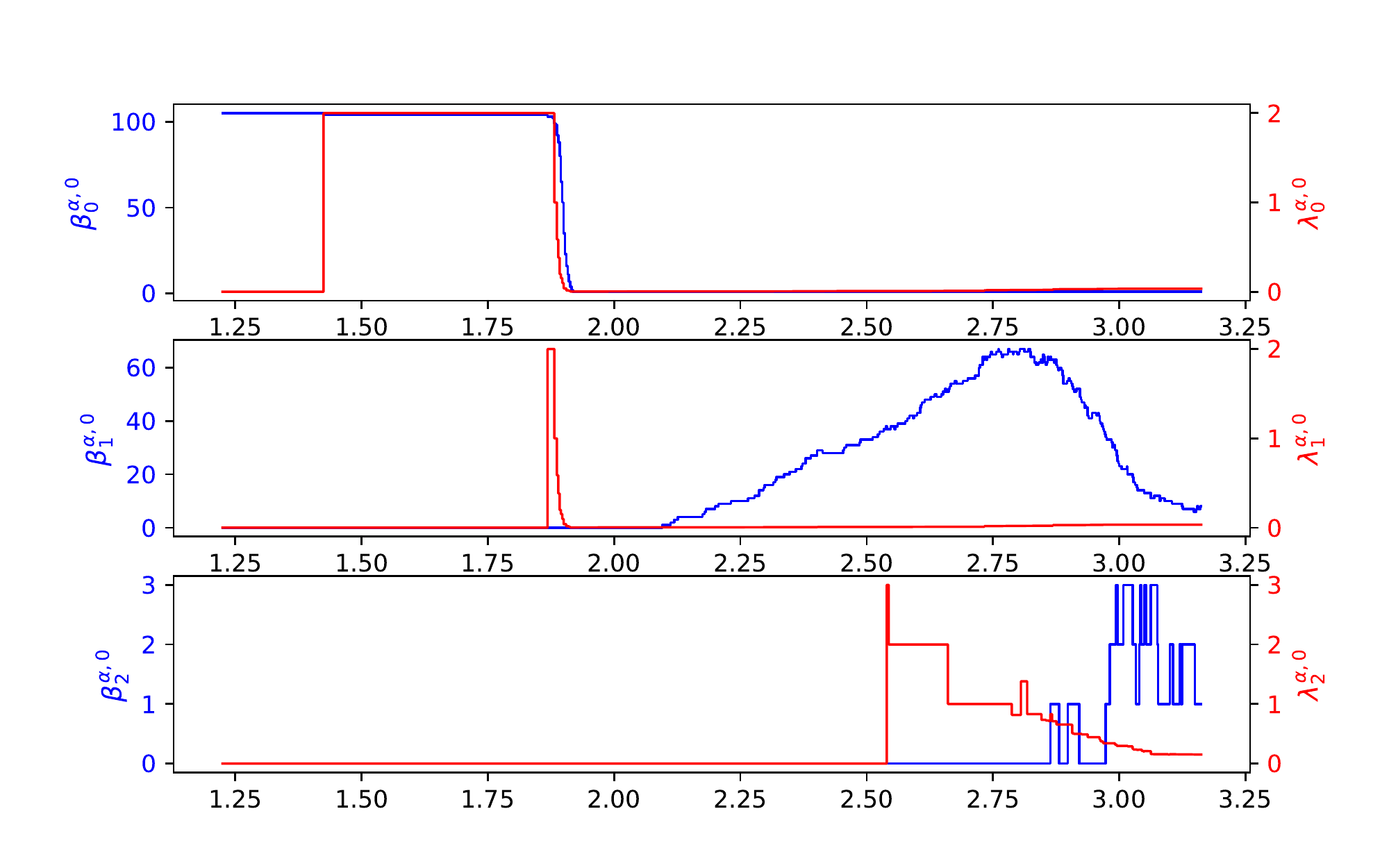}
    \caption{Illustration of the  harmonic spectra $\beta_q^{\alpha,0}$   (blue curve) and the smallest non-zero eigenvalue $\lambda_q^{\alpha,0}$ (red curve) of  PDB ID 1V70 at different filtration value $\alpha$ when $q = 0,1,2$. The $\beta_q^{\alpha,0}$ are calculated from Gudhi, DioDe, and HERMES, and $\lambda_q^{\alpha,0}$ are obtained only from HERMES. Here, the $x$-axis represents the radius filtration value $\alpha$ (unit: $\SI{}{\angstrom}$), the left-$y$-axis represents   the number of zero eigenvalues of $\mathcal{L}_q^{\alpha,0}$, and the right-$y$-axis represents   the first non-zero eigenvalue of $\mathcal{L}_q^{\alpha,0}$. Note that the harmonic spectra from three methods are indistinguishable. 
		}
    \label{fig:1V70}
\end{figure}

\begin{figure}[ht!]
    \centering
		\includegraphics[width=0.8\textwidth]{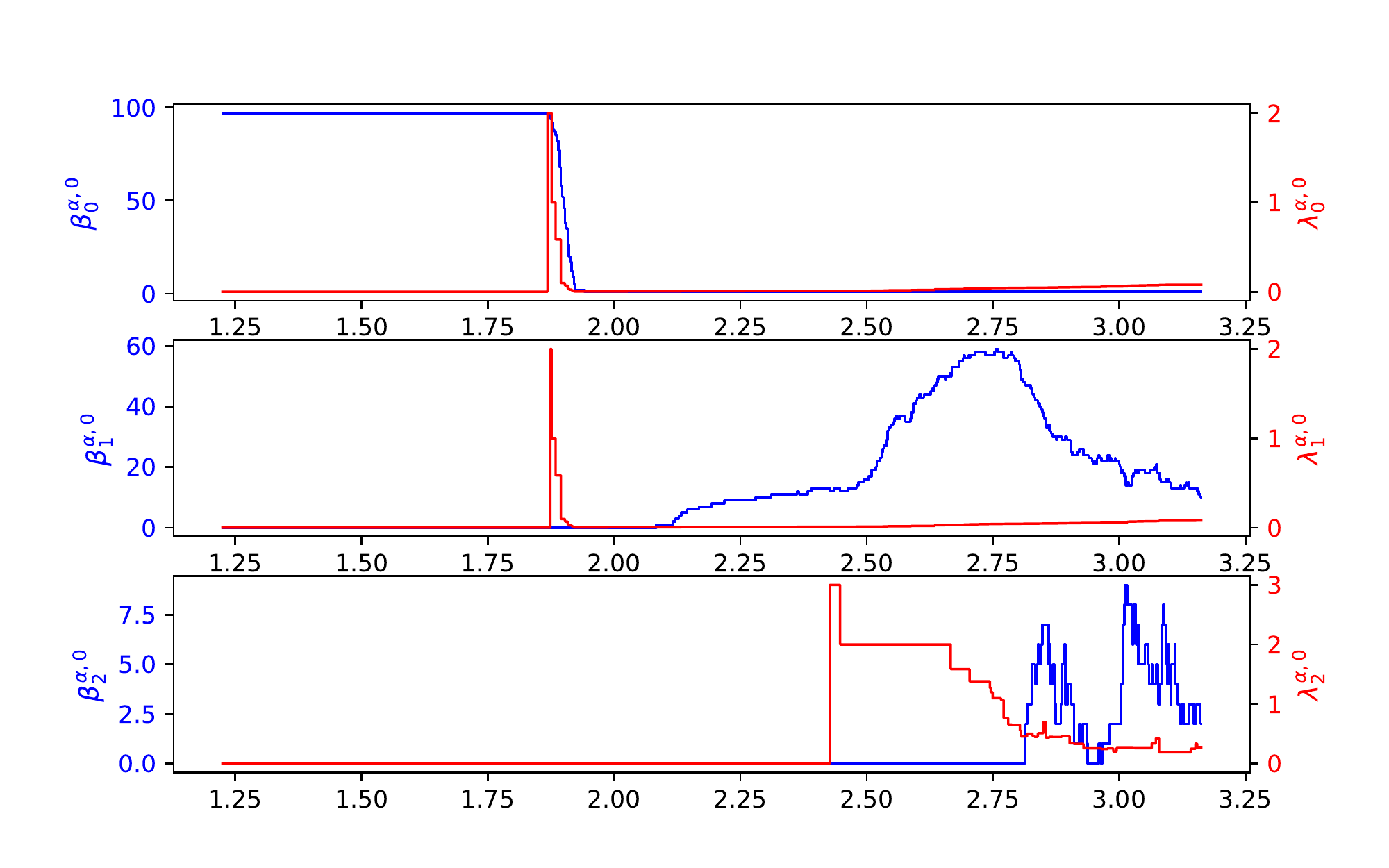}
    \caption{Illustration of the  harmonic spectra $\beta_q^{\alpha,0}$   (blue curve) and the smallest non-zero eigenvalue $\lambda_q^{\alpha,0}$ (red curve) of  PDB ID 1W2L at different filtration value $\alpha$ when $q = 0,1,2$. The $\beta_q^{\alpha,0}$ are calculated from Gudhi, DioDe, and HERMES, and $\lambda_q^{\alpha,0}$ are obtained only from HERMES. Here, the $x$-axis represents the radius filtration value $\alpha$ (unit: $\SI{}{\angstrom}$), the left-$y$-axis represents   the number of zero eigenvalues of $\mathcal{L}_q^{\alpha,0}$, and the right-$y$-axis represents   the first non-zero eigenvalue of $\mathcal{L}_q^{\alpha,0}$. Note that the harmonic spectra from three methods are indistinguishable. 
		}
    \label{fig:1W2L}
\end{figure}

\begin{figure}[ht!]
    \centering
		\includegraphics[width=0.8\textwidth]{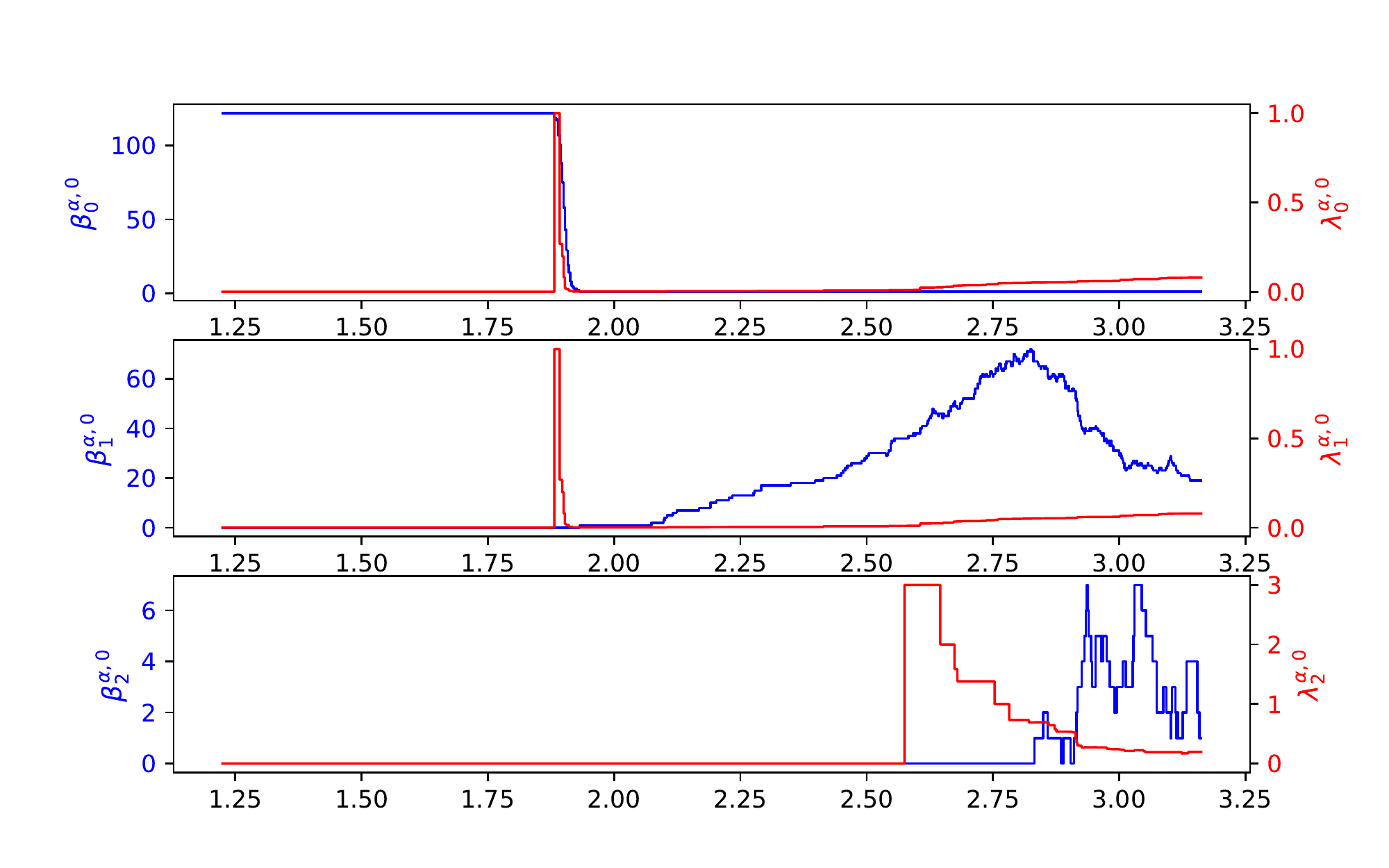}
    \caption{Illustration of the  harmonic spectra $\beta_q^{\alpha,0}$   (blue curve) and the smallest non-zero eigenvalue $\lambda_q^{\alpha,0}$ (red curve) of  PDB ID 1WHI at different filtration value $\alpha$ when $q = 0,1,2$. The $\beta_q^{\alpha,0}$ are calculated from Gudhi, DioDe, and HERMES, and $\lambda_q^{\alpha,0}$ are obtained only from HERMES. Here, the $x$-axis represents the radius filtration value $\alpha$ (unit: $\SI{}{\angstrom}$), the left-$y$-axis represents   the number of zero eigenvalues of $\mathcal{L}_q^{\alpha,0}$, and the right-$y$-axis represents   the first non-zero eigenvalue of $\mathcal{L}_q^{\alpha,0}$. Note that the harmonic spectra from three methods are indistinguishable. 
		}
    \label{fig:1WHI}
\end{figure}

\begin{figure}[ht!]
    \centering
		\includegraphics[width=0.8\textwidth]{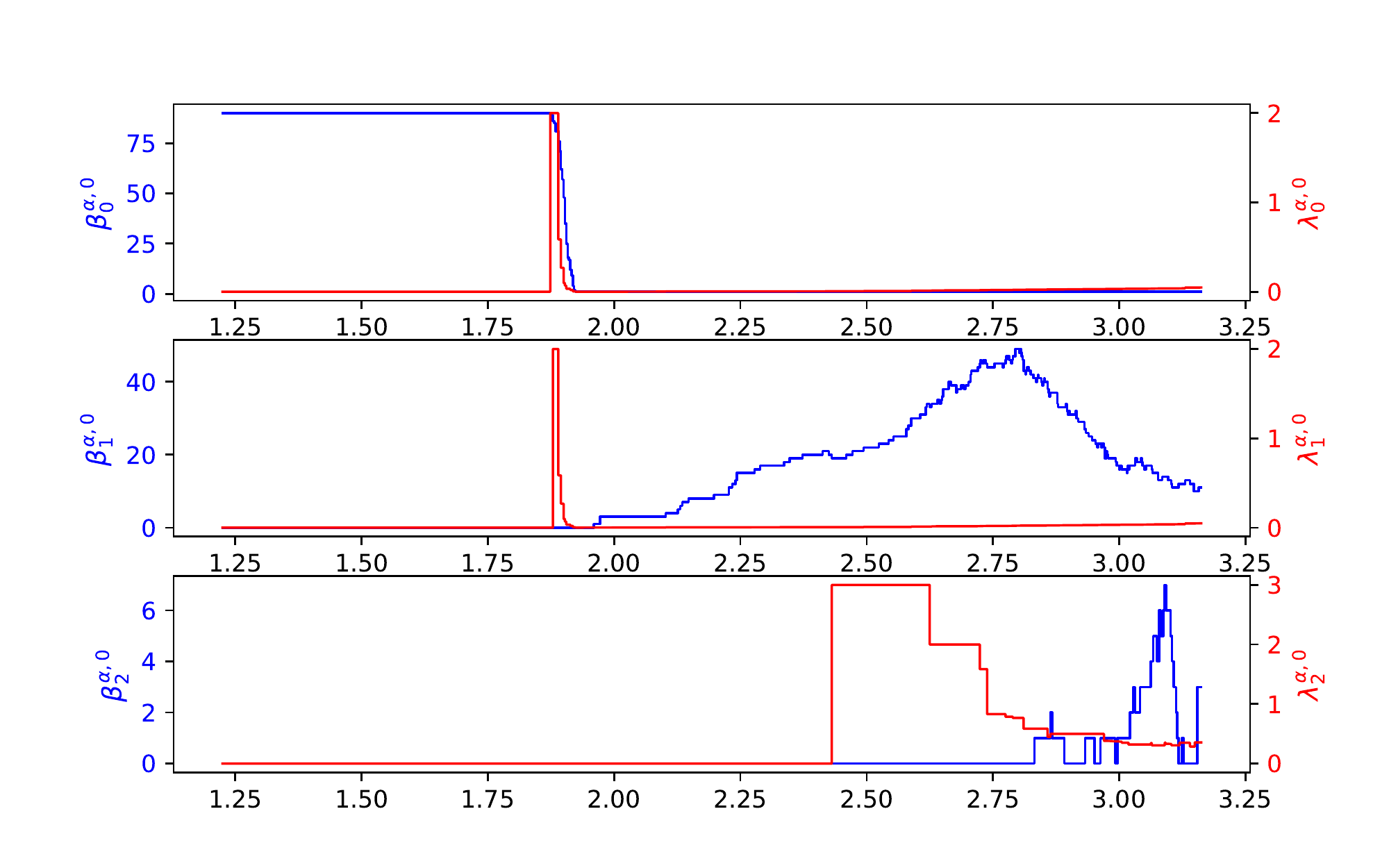}
   \caption{Illustration of the  harmonic spectra $\beta_q^{\alpha,0}$   (blue curve) and the smallest non-zero eigenvalue $\lambda_q^{\alpha,0}$ (red curve) of  PDB ID 2CG7 at different filtration value $\alpha$ when $q = 0,1,2$. The $\beta_q^{\alpha,0}$ are calculated from Gudhi, DioDe, and HERMES, and $\lambda_q^{\alpha,0}$ are obtained only from HERMES. Here, the $x$-axis represents the radius filtration value $\alpha$ (unit: $\SI{}{\angstrom}$), the left-$y$-axis represents   the number of zero eigenvalues of $\mathcal{L}_q^{\alpha,0}$, and the right-$y$-axis represents   the first non-zero eigenvalue of $\mathcal{L}_q^{\alpha,0}$. Note that the harmonic spectra from three methods are indistinguishable. 
		}
    \label{fig:2CG7}
\end{figure}

\begin{figure}[ht!]
    \centering
		\includegraphics[width=0.8\textwidth]{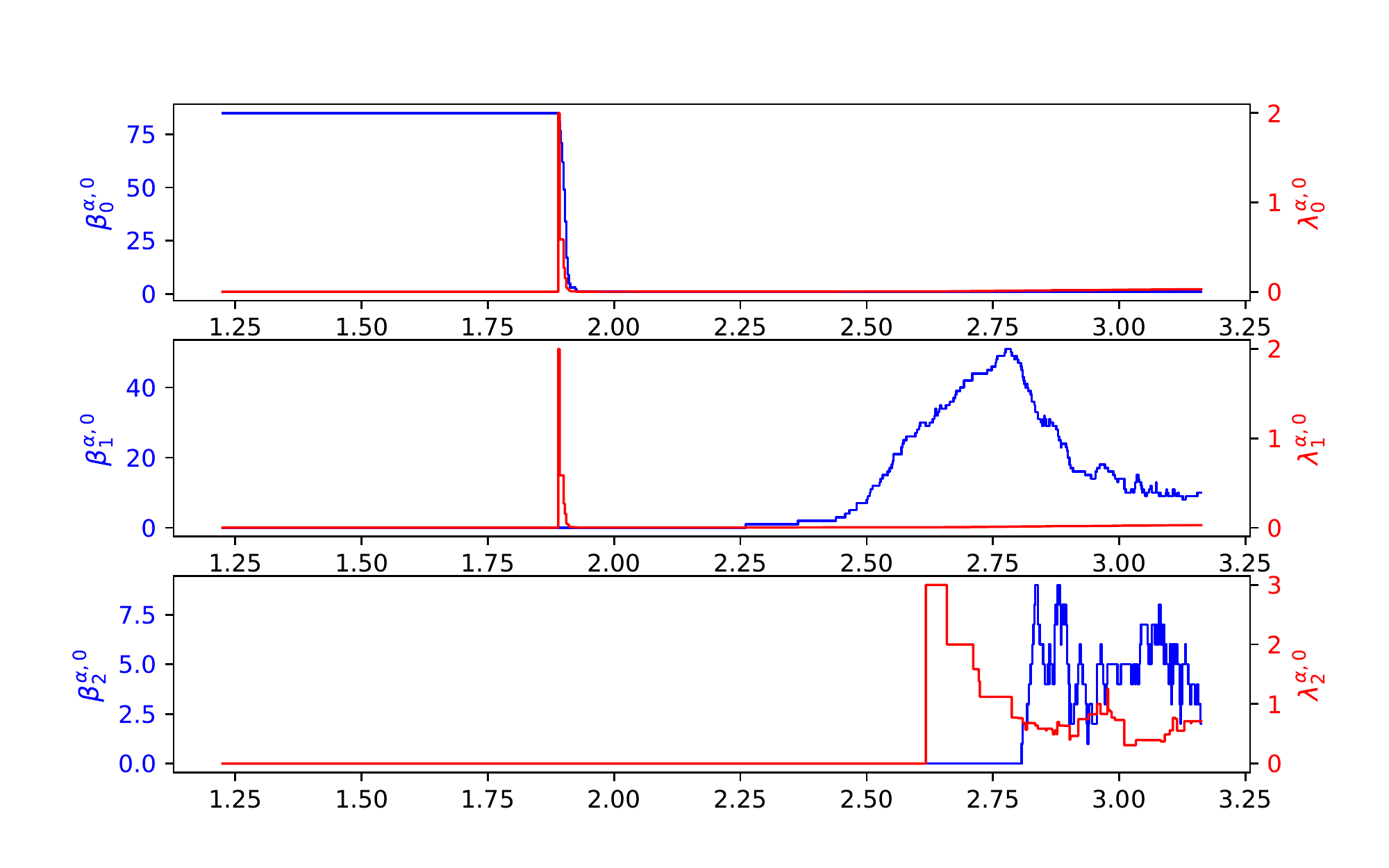}
    \caption{Illustration of the  harmonic spectra $\beta_q^{\alpha,0}$   (blue curve) and the smallest non-zero eigenvalue $\lambda_q^{\alpha,0}$ (red curve) of  PDB ID 2FQ3 at different filtration value $\alpha$ when $q = 0,1,2$. The $\beta_q^{\alpha,0}$ are calculated from Gudhi, DioDe, and HERMES, and $\lambda_q^{\alpha,0}$ are obtained only from HERMES. Here, the $x$-axis represents the radius filtration value $\alpha$ (unit: $\SI{}{\angstrom}$), the left-$y$-axis represents   the number of zero eigenvalues of $\mathcal{L}_q^{\alpha,0}$, and the right-$y$-axis represents   the first non-zero eigenvalue of $\mathcal{L}_q^{\alpha,0}$. Note that the harmonic spectra from three methods are indistinguishable. 
		}
    \label{fig:2FQ3}
\end{figure}

\begin{figure}[ht!]
    \centering
		\includegraphics[width=0.8\textwidth]{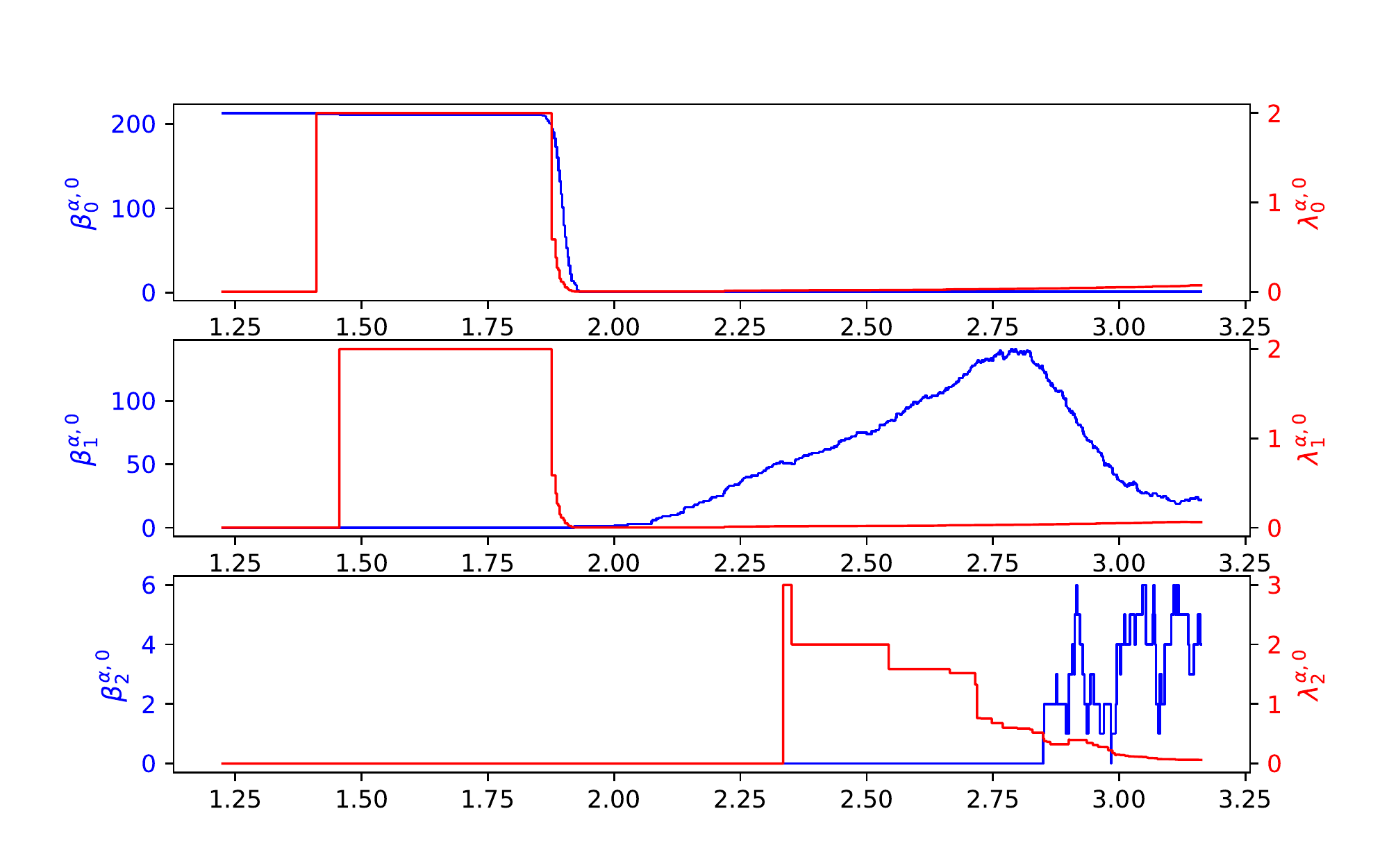}
    \caption{Illustration of the  harmonic spectra $\beta_q^{\alpha,0}$   (blue curve) and the smallest non-zero eigenvalue $\lambda_q^{\alpha,0}$ (red curve) of  PDB ID 2HQK at different filtration value $\alpha$ when $q = 0,1,2$. The $\beta_q^{\alpha,0}$ are calculated from Gudhi, DioDe, and HERMES, and $\lambda_q^{\alpha,0}$ are obtained only from HERMES. Here, the $x$-axis represents the radius filtration value $\alpha$ (unit: $\SI{}{\angstrom}$), the left-$y$-axis represents   the number of zero eigenvalues of $\mathcal{L}_q^{\alpha,0}$, and the right-$y$-axis represents   the first non-zero eigenvalue of $\mathcal{L}_q^{\alpha,0}$. Note that the harmonic spectra from three methods are indistinguishable. 
		}
    \label{fig:2HQK}
\end{figure}

\begin{figure}[ht!]
    \centering
		\includegraphics[width=0.8\textwidth]{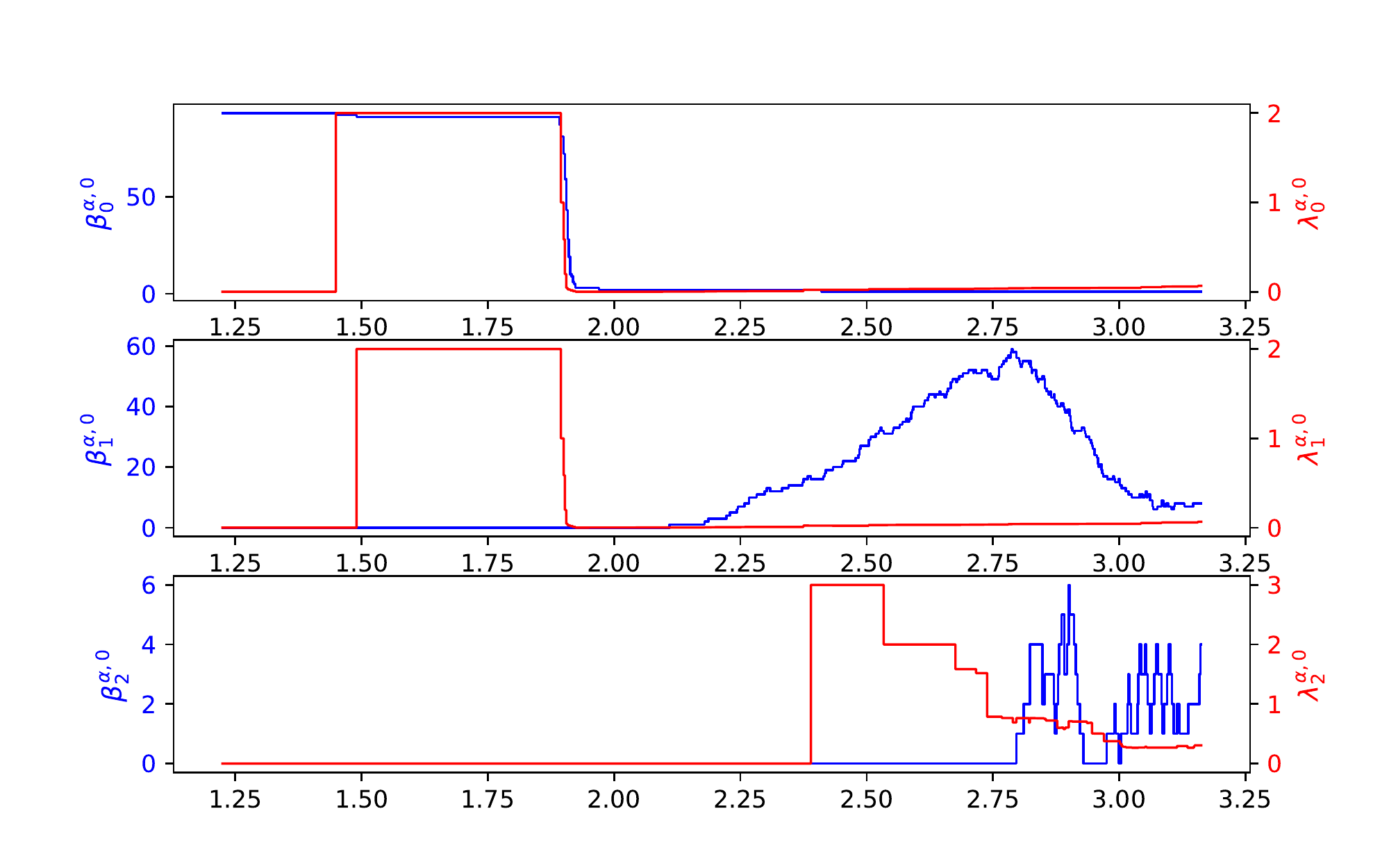}
    \caption{Illustration of the  harmonic spectra $\beta_q^{\alpha,0}$   (blue curve) and the smallest non-zero eigenvalue $\lambda_q^{\alpha,0}$ (red curve) of  PDB ID 2PKT at different filtration value $\alpha$ when $q = 0,1,2$. The $\beta_q^{\alpha,0}$ are calculated from Gudhi, DioDe, and HERMES, and $\lambda_q^{\alpha,0}$ are obtained only from HERMES. Here, the $x$-axis represents the radius filtration value $\alpha$ (unit: $\SI{}{\angstrom}$), the left-$y$-axis represents   the number of zero eigenvalues of $\mathcal{L}_q^{\alpha,0}$, and the right-$y$-axis represents   the first non-zero eigenvalue of $\mathcal{L}_q^{\alpha,0}$. Note that the harmonic spectra from three methods are indistinguishable. 
		}
    \label{fig:2PKT}
\end{figure}

\begin{figure}[ht!]
    \centering
		\includegraphics[width=0.8\textwidth]{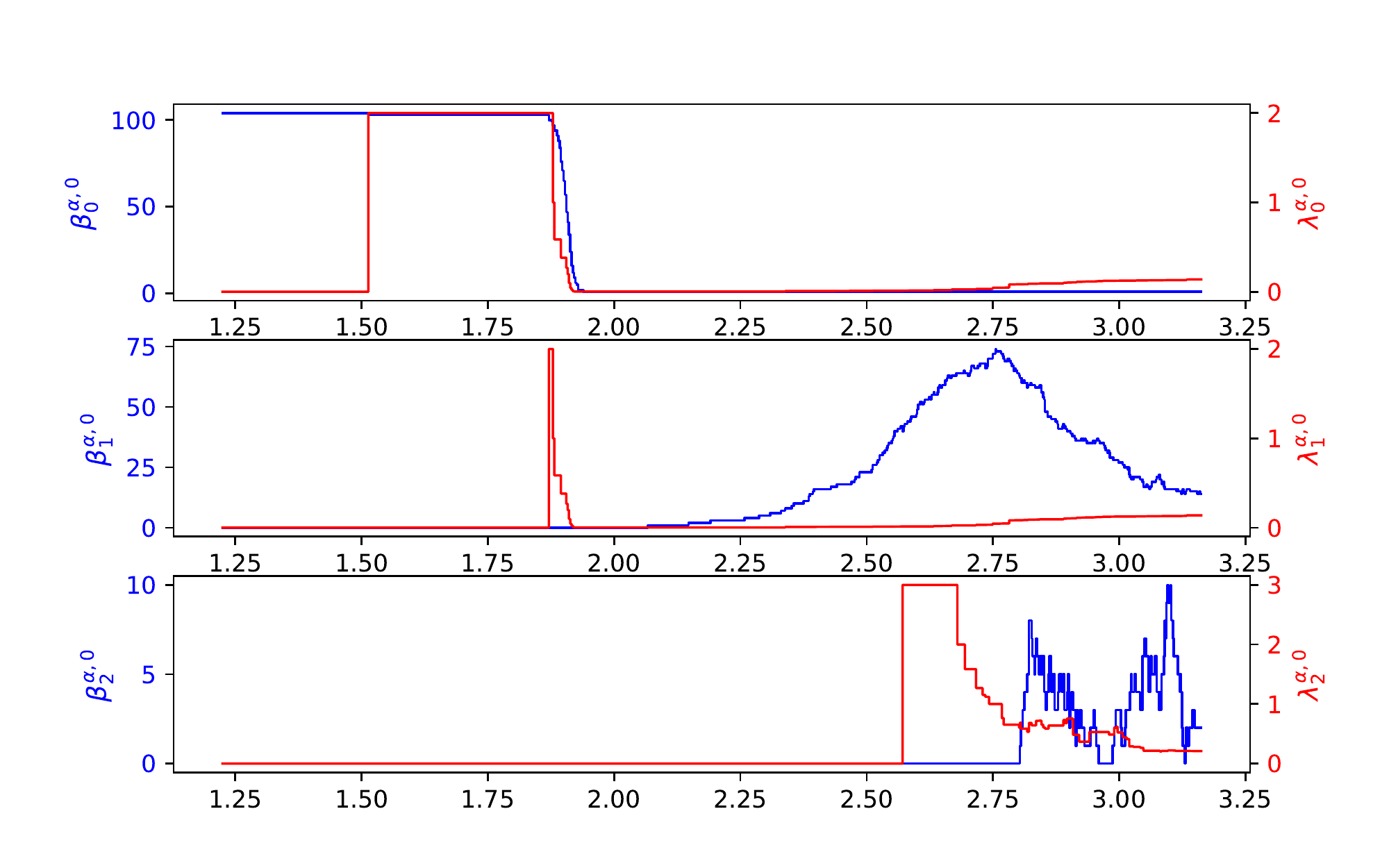}
    \caption{Illustration of the  harmonic spectra $\beta_q^{\alpha,0}$   (blue curve) and the smallest non-zero eigenvalue $\lambda_q^{\alpha,0}$ (red curve) of  PDB ID 2VIM at different filtration value $\alpha$ when $q = 0,1,2$. The $\beta_q^{\alpha,0}$ are calculated from Gudhi, DioDe, and HERMES, and $\lambda_q^{\alpha,0}$ are obtained only from HERMES. Here, the $x$-axis represents the radius filtration value $\alpha$ (unit: $\SI{}{\angstrom}$), the left-$y$-axis represents   the number of zero eigenvalues of $\mathcal{L}_q^{\alpha,0}$, and the right-$y$-axis represents   the first non-zero eigenvalue of $\mathcal{L}_q^{\alpha,0}$. Note that the harmonic spectra from three methods are indistinguishable. 
		}
    \label{fig:2VIM}
\end{figure}

\autoref{fig:1CCR} - \autoref{fig:2VIM} illustrate the harmonic spectra $\beta_q^{\alpha,0}$ ($q = 0,1,2$) of PDB IDs 1CCR, 1NKO,  1OPD, 1QTO, 1R7J, 1V70, 1W2L, 1WHI, 2CG7, 2FQ3, 2HQK, 2PKT, and 2VIM at different filtration value $\alpha$ calculated from Gudhi, DioDe, and HERMES.

\end{appendices}

\end{document}